\documentclass{QT}

\usepackage{hyperref}
\usepackage{subfig}
\usepackage[all]{xy}
\xyoption{knot}
\usepackage{rotating}

\address[gilmore@math.ucla.edu]{Allison Gilmore, Department of Mathematics, University of California Los Angeles, 520 Portola Plaza, Los Angeles, CA 90095, United States}

\numberwithin{equation}{section}

\newtheorem{theorem}{Theorem}[section]

\newtheorem{lemma}[theorem]{Lemma}

\theoremstyle{definition}
\newtheorem{definition}[theorem]{Definition}

\newtheorem{obs}[theorem]{Observation}
\newtheorem{prop}[theorem]{Proposition}

\newcommand{\noi}[0]{\noindent}
\newcommand{\isom}[0]{\cong}

\newcommand{\hfkhat}[0]{\widehat{HFK}}
\newcommand{\hfkmin}[0]{HFK^-}
\newcommand{\hfstilde}[0]{\widetilde{HFS}}

\newcommand{\cA}[0]{\mathcal{A}}
\newcommand{\cB}[0]{\mathcal{B}}
\newcommand{\cC}[0]{\mathcal{C}}
\newcommand{\cE}[0]{\mathcal{E}}
\newcommand{\cJ}[0]{\mathcal{J}}
\newcommand{\cL}[0]{\mathcal{L}}
\newcommand{\cN}[0]{\mathcal{N}}
\newcommand{\cP}[0]{\mathcal{P}}
\newcommand{\cR}[0]{\mathcal{R}}
\newcommand{\cS}[0]{\mathcal{S}}
\newcommand{\cT}[0]{\mathcal{T}}

\newcommand{\weight}[0]{\textbf{\textrm{w}}}

\title{Invariance and the knot Floer cube of resolutions}

\author[]{Allison Gilmore\thanks{The author was partially supported by NSF grant number DMS 1103801.}}

\begin{document}

\maketitle

\begin{abstract}
This paper considers the invariance of knot Floer homology in a purely algebraic setting, without reference to Heegaard diagrams, holomorphic disks, or grid diagrams. We show that (a small modification of) Ozsv\'ath and Szab\'o's cube of resolutions for knot Floer homology, which is assigned to a braid presentation with a basepoint, is invariant under braid-like Reidemeister moves II and III and under conjugation. All moves are assumed to happen away from the basepoint. We also describe the behavior of the cube of resolutions under stabilization. The techniques echo those employed to prove the invariance of HOMFLY-PT homology by Khovanov and Rozansky, and are further evidence of a close relationship between the theories. The key idea is to prove categorified versions of certain equalities satisfied by the Murakami-Ohtsuki-Yamada state model for the HOMFLY-PT polynomial.
\end{abstract}



\section{Introduction}
\label{sec:intro}

Several knot polynomials were originally categorified using a ``cube of
resolutions'' construction. Given a projection of a knot with $m$ crossings, one considers
two ways of resolving each crossing and arranges all possible resolutions
into an $m$-dimensional cube. To each vertex of the cube, one associates
a graded algebraic object (perhaps a vector space, or a module over some
commutative ring), and to each edge of the cube a map. With the correct
choices of objects and maps, the result is a chain complex whose graded
Euler characteristic is the desired knot polynomial. Khovanov's
categorification of the Jones polynomial \cite{khovjones} follows this model, employing a cube in
which the resolutions are the two possible smoothings of a
crossing. A complete resolution is then a collection of circles, to
which one associates certain vector spaces. Khovanov and Rozansky's categorification of the $\mathfrak{sl}_n$
polynomials \cite{kr1} and later the HOMFLY-PT polynomial \cite{kr2} (see also
Khovanov \cite{khovHH} and Rasmussen \cite{rasmussenonkr}) instead use a cube of resolutions built
from singularizations of crossings and oriented smoothings. The complete
resolutions in this case are a particular type of oriented planar graph. The
associated algebraic objects are modules over the ring
$\mathbb{Q}[x_0,\ldots,x_n]$, which has one indeterminate for each edge
of the graph.  In each of these theories, the chain complex was proved
to be a knot invariant by directly checking invariance under Reidemeister moves. That is,
one compares the prescribed chain complex before and after a
Reidemeister move is performed on the diagram, and constructs
a chain homotopy between the two complexes.

Knot Floer homology, which categorifies the Alexander polynomial, was
originally developed via an entirely different route. It was defined by
Ozsv\'ath and Szab\'o \cite{ozsszhfk} and by Rasmussen \cite{rasmussenthesis} as a filtration
on the chain complex of Heegaard Floer homology \cite{ozsszhf}, a
three-manifold invariant whose differentials count holomorphic disks in
the symmetric product of a surface. Knots in this theory are represented by decorating
the Heegaard diagram for a three-manifold, so invariance was proved by
checking invariance under Heegaard moves.  A second description of $\hfkhat$ was developed using grid diagrams in \cite{SWcomb} and \cite{MOScomb}. This definition is fully combinatorial and its invariance was proved combinatorially in \cite{MOSTcomb} by checking grid moves.

In 2007, Ozsv\'ath, Szab\'o and Stipsicz \cite{ozsszstipsing} described a version of knot Floer homology for singular knots that is related to the theory for classical knots by a skein exact sequence. In general, knot Floer homology for singular knots involves holomorphic disk counts, but it can be made combinatorial with a suitable choice of twisted coefficients and a particular Heegaard diagram.  Using this version of the theory for singular knots and iterating the skein exact sequence allowed Ozsv\'ath and Szab\'o \cite{ozsszcube} to calculate knot Floer homology using a cube of resolutions. Their construction is, in the end, fully algebraic. Compared to the grid diagram formulation of knot Floer homology, it has the advantage of being not only combinatorial, but conceptually grounded.

This paper addresses the question of invariance for knot Floer homology within the algebraic setting of the cube of resolutions. Our invariance result is weaker than might be expected---the Reidemeister I move is missing and other moves must avoid a basepoint---but our methods advance the project of understanding knot Floer homology from an algebraic perspective.  We do not rely on Heegaard diagrams, holomorphic disks, or any of the usual geometric input. We also make no reference to grid diagrams. The main idea is to prove categorified versions of certain equalities satisfied by the Murakami-Ohtsuki-Yamada model of the HOMFLY-PT polynomial. 

A thorough understanding of knot Floer homology from the algebraic perspective should have several applications. Perhaps most immediately, it should clarify the relationship between knot Floer homology and other knot homologies. It also suggests a relationship between the knot Floer cube of resolutions and the various constructions of the Alexander polynomial via representation theory~\cite{kauffmansaleur},\cite{rozanskysaleur},\cite{murakamialex},\cite{murakamialex2},\cite{deguchiakutsu},\cite{viroquantumrel}. A full categorification of such constructions would extend knot Floer homology to a tangle invariant.

We use a small modification of the cube of resolutions construction described by Ozsv\'ath and Szab\'o \cite{ozsszcube}, which still produces chain complexes with homology $\hfkhat$ and $\hfkmin$. 
We begin with a projection of a knot as a closed braid, which we decorate with a number of extra bivalent vertices. We also place a basepoint on one of the outermost edges. (Equivalently, we cut the outermost strand to obtain a $(1,1)$-tangle in braid position.) The result is a \emph{layered braid diagram} $D$. 
We form a cube of resolutions by singularizing or smoothing each crossing of the projection. We then assign a graded algebra $\cA_I(D)$ to each resolution and arrange these into a chain complex $C(D)$. These objects are defined precisely in (\ref{corneralg}) and (\ref{cubedfn}) of Section~\ref{sec:definitions}.  Our main result is an algebraic proof of 
\begin{theorem}
\label{thm:invc}
Let $D$ be a layered braid diagram for a knot. The chain complex $C(D)$, up to chain homotopy equivalence and twisting by certain endomorphisms of the ground ring, is invariant under braid-like Reidemeister moves II and III and under conjugation. All moves are assumed to avoid the basepoint. 
 \end{theorem}
 \noi Theorem~\ref{thm:invc} holds with coefficients in $\mathbb{Z}$. It is stated in full detail in Section~\ref{sec:invc}. Changing to $\mathbb{Z}_2$ coefficients, we identify $H_\ast(C(D))$ with $\hfkmin$ and a reduced version of $C(D)$ with $\hfkhat$ in Proposition~\ref{prop:computeshfk}. We expect that $H_\ast(C(D))$ in fact computes knot Floer homology with integer coefficients, but do not pursue this point here.

The proof of Theorem~\ref{thm:invc} is modeled on Khovanov and Rozansky's invariance proof for HOMFLY-PT homology in \cite{kr2}. Specifically, we prove categorified versions of the braid-like Murakami-Ohtsuki-Yamada (MOY) relations shown in Figure~\ref{fig:braidmoy}. Murakami, Ohtsuki, and Yamada~\cite{moy} present a state sum model for the $\mathfrak{sl}(n)$ polynomials, which extends to a model for the HOMFLY-PT polynomial. To compute the polynomial invariant of a knot, one sets up a weighted sum of combinatorially defined states on certain graphs obtained from a projection of the knot. (Singular knots are equivalent to a subclass of these graphs.) The MOY model is purely combinatorial, but is derived from the Reshetikhin-Turaev~\cite{rt} recipe for producing polynomial invariants from the representation theory of quantum groups. Lemma~\ref{lemma:reid2} and Lemma~\ref{reid3splitting101} both categorify the top line of Figure~\ref{fig:braidmoy}, with the $q$ and $q^{-1}$ to be interpreted as shifts in the Alexander grading on $\cA$. Lemma~\ref{reid3splitting111} categorifies the middle and bottom lines of Figure~\ref{fig:braidmoy}.

\begin{figure}[htbp]
\begin{center}
\input{braidmoy}
\caption{Murakami-Ohtsuki-Yamada relations used in the proof of Theorem~\ref{thm:invc}.}
\label{fig:braidmoy}
\end{center}
\end{figure}

Our use of diagrams with basepoints and basepoint-avoiding moves is unsurprising from the representation theoretic point of view. Although the HOMFLY-PT polynomial specializes to the Alexander polynomial, the MOY model does not specialize to a model of the Alexander polynomial. Specializing the model so that it satisfies the Alexander polynomial's skein relation produces the zero polynomial. This unfortunate fact is not a defect of the MOY model, but a fundamental feature of the representation theory underlying the polynomial knot invariants. The Reshetikhin-Turaev method for constructing knot polynomials produces the zero polynomial when applied to the quantum groups related to the Alexander polynomial ($\mathcal{U}_q(\mathfrak{gl}(1\vert 1))$ or $\mathcal{U}_q(\mathfrak{sl}(2))$ at a root of unity).\footnote{These observations go back at least to Reshetikhin~\cite{reshetikhinalex} in the early 1990s. The problem appears to be related to the fact that the Alexander polynomial vanishes on split links.} The standard means of circumventing this dilemma has been to cut a strand and work with $(1,1)$-tangles up to isotopies fixing the endpoints~\cite{kauffmansaleur},\cite{rozanskysaleur},\cite{murakamialex},\cite{murakamialex2},\cite{deguchiakutsu},\cite{viroquantumrel}, which correspond exactly to knots with basepoints up to isotopies fixing the basepoint. Of course, the location of the cut / basepoint should not matter for a construction of knot Floer homology, but the behavior of $C(D)$ under movement of the basepoint appears to be complicated. We hope to return to this issue in future work.

We describe the behavior of $C(D)$ under stabilization in Proposition~\ref{reid1prop}. The description is separate from Theorem~\ref{thm:invc} because it requires leaving the layered braid setting. Unfortunately, the categorified MOY relations that underlie our arguments for Reidemeister moves II and III do not appear to extend to this more general setting. A new approach will probably be needed to prove that $C(D)$ (or an appropriate generalization thereof) is invariant under stabilization.

One final note is in order regarding the limitations of the cube of resolutions construction. We are not aware of any intrinsic justification for the use of braid diagrams. They appear here for the same reason as in Ozsv\'ath and Szab\'o's original construction: the holomorphic disk counting required to link the cube of resolutions to knot Floer homology is tractable only for braid diagrams. Interestingly, HOMFLY-PT homology also requires braid presentations.

This paper is organized as follows. Section~\ref{sec:definitions} describes the modified construction of the cube of resolutions needed to incorporate layered braid diagrams.  Section~\ref{sec:nonlocalrelations} examines in detail the non-local relations involved in the definition of the algebra associated to a resolution. These relations are a key difference between the cube of resolutions theories for $HFK$ and for HOMFLY-PT homology. Section~\ref{sec:removemark} establishes a technical proposition allowing us to remove sets of bivalent vertices under certain conditions. The next sections address Reidemeister II, Reidemeister III, conjugation, and stabilization in turn. Section~\ref{sec:computeshfk} verifies that the cube of resolutions defined here computes knot Floer homology.

\subsection*{Acknowledgments}
I would like to thank Peter Ozsv\'ath for suggesting this problem, for patiently explaining the material in \cite{ozsszcube}, and for always being a supportive and encouraging advisor. I am also grateful to Mikhail Khovanov for several thought-provoking conversations.  I thank Daniel Krasner for illuminating the technical material in \cite{kr2} and John Baldwin for helpful conversations about an earlier version of this paper. Finally, I am grateful to an anonymous referee for a very thorough review of that prior version, and to Robert Lipshitz and Ciprian Manolescu for their encouragement during revisions.

\section{Definitions: Cube of Resolutions for HFK}
\label{sec:definitions}

\begin{figure}[htbp]
\begin{center}
\input{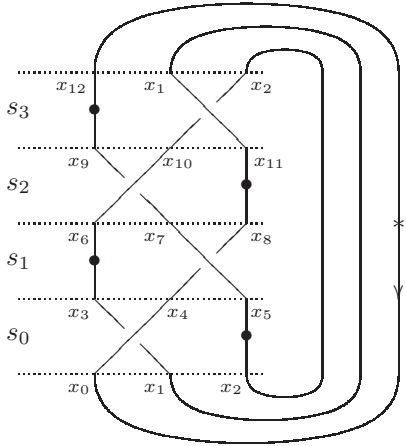}
\caption{A layered braid diagram for the figure 8 knot.}
\label{notation}
\end{center}
\end{figure}

We begin with an oriented braid-form projection $D$ of an oriented knot $K$ in $S^3$. Let $b$ refer to the number of strands in $D$ (which is not necessarily the braid index of $K$). Subdivide one of the outer edges of $D$ by a basepoint~$\ast$. Isotoping $D$ as necessary, fix an ordering on its crossings so that $D$ is the closure of a braid diagram that is a stack of $m+1$ horizontal layers, each containing a single crossing and $b-2$ vertical strands. Label the horizontal layers $s_0,\ldots,s_m$. This amounts to choosing a braid word for $D$. In each horizontal layer, add a bivalent vertex to each strand that is not part of the crossing. Finally, label the edges of $D$ by $0,\ldots,n$, where $n=(m+1)b$, such that 0 is the edge coming out from the basepoint  and $n$ is the edge pointing into the basepoint. A braid diagram in this form will be called a \emph{layered braid diagram} for $K$. See Figure~\ref{notation} for an example of a layered braid diagram of the figure 8 knot.  Although Ozsv\'ath and Szab\'o \cite{ozsszcube} use closed braid diagrams with basepoints in their definition of the knot Floer cube of resolutions, they do not require diagrams to be layered. This refinement appears to be critical to our proof of Proposition~\ref{prop:removemark} and necessary for the proofs of the categorified MOY relations (Lemmas~\ref{lemma:reid2},~\ref{reid3splitting101},~\ref{reid3splitting111}) that underlie Reidemeister moves II and III.

\begin{figure}[htbp]
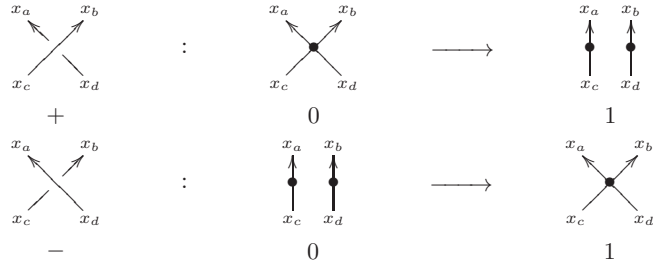

\begin{center}
{\scalebox{.9}{\input{positiveresolutionlabels.tex}}}
{\scalebox{.9}{\input{negativeresolutionlabels.tex}}}
\caption{Notation for the singularization and smoothing of a positive (respectively negative) crossing.}
\label{resolutionlabels}
\end{center}
\end{figure}

Each crossing in a knot projection can be singularized or smoothed.  To singularize the crossing in layer $s_i$, replace it by a 4-valent vertex and retain all edge labels.  To smooth the crossing in layer $s_i$, replace it with two vertical strands with one bivalent vertex on each, and retain all edge labels. Figure~$\ref{resolutionlabels}$ illustrates these labeling conventions.  

A resolution of a knot projection is a diagram in which each crossing has been singularized or smoothed.  Alternatively, it is a planar graph in which each vertex is either (1) 4-valent with orientations as in Figure~\ref{resolutionlabels}, or (2) bivalent with one incident edge oriented towards the vertex and the other oriented away from the vertex.  For a positive crossing, declare the singularization to be the 0-resolution and the smoothing to be the 1-resolution.  For a negative crossing, reverse these labels. A resolution of a knot projection can then be specified by a multi-index of 0s and 1s, generically denoted $\epsilon_0\ldots\epsilon_m$, or simply $I$, which we will think of as a vertex of a hypercube.  Considering all possible singularizations and smoothings of all crossings, we obtain a cube of resolutions for the original knot projection. The homological grading on the cube will be given by collapsing diagonally; that is, by summing $\epsilon_0+\cdots+\epsilon_m$.

Let $\cR=\mathbb{Z}[t^{-1},t]$ and $\underline{x}(D)$ denote a set of formal variables $x_0,\ldots,x_n$ corresponding to the edges of $D$. Define the \emph{edge ring} of $D$ to be $\cR[\underline{x}(D)]$, which we will abbreviate to $\cR[\underline{x}]$ if $D$ is clear from context.   To each vertex of the cube of resolutions, we will associate an $\cR$-algebra $\cA_I(D)$, which is a quotient of the edge ring by an ideal defined by combinatorial data in the $I$-resolution of $D$. To each edge of the cube, we will associate a map.  Together with proper choices of gradings, these data define a chain complex of graded algebras over $\cR[\underline{x}(D)]$.  We will sometimes need to complete $\cR$ or $\cR[\underline{x}(D)]$ with respect to $t$, meaning that we will allow Laurent series in $t$ with coefficients in $\mathbb{Z}$ or $\mathbb{Z}[\underline{x}(D)]$, respectively. Denote these completions $\widehat{\cR}$ and $\widehat{\cR[\underline{x}(D)]}$, respectively.  More specifically, the description of the cube of resolutions' behavior under stabilization requires extending the ground ring to $\widehat{\cR}[\underline{x}(D)]$ and the identification of the homology of $C(D)$ with knot Floer homology requires extending to $\widehat{\cR[\underline{x}(D)]}$ (as well as passing to $\mathbb{Z}_2$ coefficients).

\subsection{Algebra associated to a resolution}
\label{sec:algdfn}
The algebra associated to the $I$-resolution of the knot projection $D$, which we will denote $\cA_I(D)$, is the quotient of the edge ring by the ideal generated by the following three types of relations.

\begin{enumerate}
\item Linear relations associated to each vertex. 
\begin{align*}
t(x_a+x_b)&-(x_c+x_d)\quad\text{to}\quad{\xy
(-4,-4)*{}="bl";
(4,-4)*{}="br";
(-4,4)*{}="tl";
(4,4)*{}="tr";
(0,0)*{\bullet}="mm";
{\ar "bl";"tr"};
{\ar "br";"tl"};
(-5,5.5)*{x_a};
(5,5.5)*{x_b};
(-5,-5.5)*{x_c};
(5,-5.5)*{x_d};
\endxy}\\
\\
tx_{a}&-x_c\quad\text{to}\quad
{\xy
(0,-6)*{}="b";
(0,6)*{}="t";
(0,0)*{\bullet}="m";
{\ar "b";"t"};
(2.5,0)*{};
(3,5)*{x_{a}};
(2.5,-5)*{x_c};
\endxy}\\
\end{align*}
\item Quadratic relations associated to each 4-valent vertex.  \begin{align*}
t^2x_ax_b-x_cx_d\quad\text{to}\quad
{\xy
(-4,-4)*{}="bl";
(4,-4)*{}="br";
(-4,4)*{}="tl";
(4,4)*{}="tr";
(0,0)*{\bullet}="mm";
{\ar "bl";"tr"};
{\ar "br";"tl"};
(-5,5.5)*{x_a};
(5,5.5)*{x_b};
(-5,-5.5)*{x_c};
(5,-5.5)*{x_d};
\endxy}
\end{align*}  Note that this relation can always be rewritten in four different ways by combining with the linear relation corresponding to the same vertex: \begin{eqnarray*}(tx_a-x_c)(x_d-tx_a)\quad&\text{or}&\quad(tx_b-x_c)(x_d-tx_b)\quad\text{or}\quad\\(tx_a-x_c)(tx_b-x_c)\quad&\text{or}&\quad(tx_a-x_d)(tx_b-x_d).\end{eqnarray*}
\item Non-local relations associated to sets of vertices in the resolved diagram. These have several equivalent definitions, which will be explored in detail in Section~\ref{sec:nonlocalrelations}. Denote the ideal generated by non-local relations in $I$-resolution of $D$ by $\mathcal{N}_I(D)$ or simply $\cN_I$.
\end{enumerate}

We refer to the linear and quadratic relations as the local relations. Let $\cL$ denote the ideal they generate together, and $\cL_i$ denote the ideal generated by the local relations in layer $s_i$.  Then the algebras that belong at the corners of the cube of resolutions are
\begin{equation}
\label{corneralg}
\cA_I(D)=\frac{\cR[x_0,\ldots,x_n]}{\cL+\cN_I}.
\end{equation}

Throughout this paper, we will use ``$\equiv$'' to indicate that two polynomials in the edge ring are equivalent up to multiplication by units in $\cR[\underline{x}(D)]/\cL$. Such polynomials are equivalent in the sense that they generate the same ideal in $\cR[\underline{x}(D)]/\cL$. We will represent generating sets for ideals as single-column matrices.  The entries of the matrices are elements of the edge ring. The matrices can be manipulated using row operations without changing the ideal they generate because the ideal $(a,b)$ is identical to the ideal $(a,b+sa)$ for any unit $s\in\mathcal{R}$ and $a,b\in \cR[\underline{x}(D)]$. Also, when we see a row of the form $a-b$ in a matrix, we can replace $b$ by $a$ in all other rows and eliminate $b$ from the edge ring. This will not change the quotient of the edge ring by the ideal of relations.  Although the matrix manipulations in the following sections look very similar to those in \cite{kr2} and \cite{rasmussenonkr}, the matrices here do not formally represent matrix factorizations.

Let $S$ denote the $I$-resolution of $D$, treated as a singular knot. The algebra $\mathcal{A}_I(D)$ is a twisted version of the singular knot Floer homology of $S$. More precisely, there is a chain complex $C^\prime(S)$ over $\mathbb{Z}_2[\underline{x}(D)]$ mentioned in~\cite[Section 4]{ozsszstipsing} that is defined using Heegaard diagrams for singular knots and a differential that counts holomorphic disks. With appropriately twisted coefficients, the homology of $C^\prime(S)$ with respect to the holomorphic disk counting differential is $\cA_I(D)\otimes_{\mathbb{Z}}\mathbb{Z}_2$~\cite[Theorem 3.1]{ozsszcube}. The complex $C^\prime(S)$ is a generalization of those that are the focus of~\cite{ozsszstipsing}. Setting $x_a=x_b$ at each 4-valent vertex and $x_0$ to zero in $C^\prime(S)$, then taking homology with respect to the same differential, gives the theory called $HFS$ in~\cite{ozsszstipsing} while setting all of the edge variables to zero before taking homology gives the theory they call $\hfstilde$. It is proved in~\cite[Theorem 1.3, Section 5]{ozsszstipsing} that $HFS$ is completely determined by its Euler characteristic, while $\hfstilde$ contains additional information.

\subsection{Differential}
\label{sec:diffl}

An edge of the cube of resolutions goes between two resolutions that differ at exactly one crossing. To an edge that changes the $i^{th}$ crossing, we associate a map $\mathcal{A}_{\epsilon_0\ldots 0 \ldots\epsilon_m}(D)\longrightarrow\mathcal{A}_{\epsilon_0\ldots 1\ldots\epsilon_m}(D)$.  If $s_i$ was positive in the original knot projection, then the edge goes from a diagram containing the singularization of $s_i$ to a diagram containing its smoothing. The ideal of relations associated to the singularized crossing is contained in the ideal of relations associated to the resolved crossing (see Observation~\ref{singsmoothobs}), so $\mathcal{A}_{\epsilon_0\ldots 1\ldots\epsilon_m}(D)$ is a quotient of $\mathcal{A}_{\epsilon_0\ldots 0 \ldots\epsilon_m}(D)$. The corresponding map in this case will be the quotient map. If $s_i$ was negative in the original knot projection, then the edge goes from the smoothing to the singularization of $s_i$. The corresponding map in this case will be multiplication by $tx_a-x_d$, or equivalently by $tx_b-x_c$, where the crossing $s_i$ is labeled as in Figure~\ref{resolutionlabels}.

We assign signs to the edge maps according to the following procedure. First, for each edge we have a choice of two maps, $tx_a-x_d$ or $tx_b-x_c$, which are equal up to a sign. We assign a negative sign to those of the form $tx_a-x_d$ and a positive sign to those of the form $tx_b-x_c$. This produces a cube in which all squares commute. We then change signs on certain edges to turn commuting squares into anti-commuting squares. Let $\epsilon_0,\ldots,\epsilon_m\in\{0,1\}$, $\epsilon_i=0$, and $\hat{\epsilon}_i=1$. Then the arrow from vertex $\epsilon_0\ldots\epsilon_i\ldots\epsilon_m$ to $\epsilon_0\ldots\hat{\epsilon}_i\ldots\epsilon_m$ keeps its existing sign if $\epsilon_0+\cdots+\epsilon_{i-1}\equiv 0\mod 2$ and changes sign otherwise.

We have now assembled all of the pieces needed to define the chain complex $(C(D),d)$ referred to in Theorem~\ref{thm:invc}.  Let
\begin{equation}
\label{cubedfn}
C(D)=\bigoplus_{I\in\{0,1\}^{m+1}}\cA_I(D)
\end{equation}
with total differential $d$ the sum of all edge maps and homological grading given by $\epsilon_0+\cdots+\epsilon_m$. This is the chain complex that computes $\hfkmin$ (see Proposition~\ref{prop:computeshfk}). There is also a reduced version of this chain complex obtained by setting $x_0$ to zero in each $\cA_I(D)$. Its homology computes $\hfkhat$.

\subsection{Gradings} The chain complex $C(D)$ comes equipped with an additional grading called the Alexander grading. Let $\mathcal{R}$ be in grading 0 and each edge variable $x_i$ in grading -1.  The relations used to form $\mathcal{A}_I(D)$ are homogeneous with respect to this grading, so it descends from the edge ring to a grading on $\mathcal{A}_I(D)$ (called $A_0$ in \cite{ozsszcube}). To symmetrize, adjust upwards by a factor of $\frac{1}{2}\left(m_\times-b+1\right)$, where $m_\times$ is the number of singular points in the $I$-resolution of $D$ and $b$ is the number of strands in $D$. Call this the internal grading, $A_I$, on $\cA_I(D)$.

The Alexander grading on $\mathcal{A}_I(D)$ as a summand of the cube $C(D)$ is further adjusted from the internal grading by
\begin{equation*}
\label{cubealexgrading}
A=A_I+\frac{1}{2}\left(-m_-+\sum_{i=0}^m\epsilon_i\right),
\end{equation*}
where $\epsilon_0,\ldots,\epsilon_m$ are the components of the multi-index $I$ and $m_-$ is the number of negative crossings in $D$.  This grading $A$ is the final Alexander grading on the complex $C(D)$.

\subsection{Invariance}
\label{sec:invc}
With these definitions in place, we may now state Theorem~\ref{thm:invc} precisely.
\begin{theorem}
Let $D$ and $D^\prime$ be layered braid diagrams that are related by a finite sequence of braid-like Reidemeister moves II and III and conjugation. Assume that all such moves avoid the basepoint. Let $x_0$ be the initial edge of $D$ and $D^\prime$. Let $\psi_i$ be the endomorphism of $\cR$ that takes 1 to 1 and $t$ to $t^i$.

Consider $C(D)$ and $C(D^\prime)$ as complexes of graded $\cR[x_0]$-modules. There is a complex $C(D,D^\prime)$ of graded $\cR[x_0]$-modules such that $C(D)$ and $C(D^\prime)$ are each chain homotopy equivalent to twistings of $C(D,D^\prime)$ by a composition of finitely many $\psi_i$. 

Alternatively, there is a collection of positive integers $m_1,\ldots,m_N$ such that $C(D)\otimes_\cR\cR[t^{1/m_1},\ldots,t^{1/m_N}]$ and $C(D^\prime)\otimes_\cR\cR[t^{1/m_1},\ldots,t^{1/m_N}]$, considered as complexes of graded $\cR[t^{1/m_1},\ldots,t^{1/m_N}][x_0]$-modules, are related by chain homotopy equivalence and twisting by various $\psi_i$ and their inverses.
\end{theorem}

\section{Non-local Relations}
\label{sec:nonlocalrelations}

\begin{figure}[htbp]
\begin{center}
\input{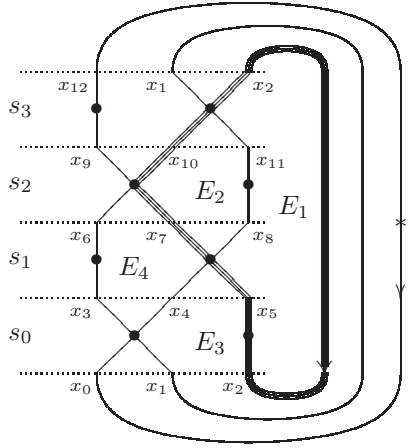}
\caption{Singularization of the minimal braid presentation of the figure 8 knot with edges labeled $x_0,\ldots,x_{12}$ and orientations consistent with those in Figure~\ref{notation}. The bold line shows a cycle whose corresponding non-local relation is $t^8x_1x_9-x_4x_6$. Elementary regions are labeled $E_1,\ldots,E_4$.  The coherent region $E_1\cup E_2$ produces the same non-local relation as the cycle in bold, as does the subset consisting of the bivalent vertex in $s_0$, the 4-valent vertex in $s_1$, all vertices in $s_2$, and the 4-valent vertex in $s_3$.}
\label{fig:nonlocal}
\end{center}
\end{figure}

We collect here three equivalent definitions of the non-local relations used in the description of the algebra $\cA_I(D)$, along with several straightforward observations that will nonetheless be very useful in later arguments. Figure~\ref{fig:nonlocal} will serve as a source of examples throughout.

First, we may generate $\mathcal{N}_I$ by associating a relation to each cycle (simple closed path) in the resolved diagram that does not pass through the basepoint and that is oriented consistently with $D$.  

\begin{definition}[Cycles]
\label{dfncycles}
Let $Z$ be a simple closed path in the $I$-resolution of $D$ that does not pass through the basepoint and is oriented consistently with $D$.  Let $R_Z$ be the region it bounds in the plane, containing the braid axis.  The \emph{weight} $\emph{\weight}(Z)$ of $Z$ is twice the number of 4-valent vertices plus the number of bivalent vertices in the closure of $R_Z$.  The \emph{non-local relation associated to $Z$} is $$t^{\emph{\weight}(Z)}w_{\text{out}}-w_{\text{in}},$$ where $w_{\text{out}}$ (respectively $w_{\text{in}}$) is the product of the edges incident to exactly one vertex of $Z$ that lie outside of $R_Z$ and that point out of (respectively into) the region.
\end{definition}

\noi Figure~\ref{fig:nonlocal} shows a cycle in the singularized figure 8 knot with associated relation $t^8x_1x_9-x_4x_6$.

A slightly different definition derives a generating set for $\mathcal{N}_I$ from certain regions in the complement of the $I$-resolution of $D$. First define the \emph{elementary regions} in the $I$-resolution of $D$ to be the connected components of its complement in the plane, except for the two components that are adjacent to the basepoint. For example, there are four elementary regions in the singularized figure 8 shown in Figure~\ref{fig:nonlocal}.

Since $D$ is assumed to be in braid position, the elementary regions can be partially ordered based on which two strands of $D$ they lie between. Label the strands of $D$ from 1 (innermost, nearest the braid axis) to $b$ (outermost, nearest the non-compact region).  Then $E_i<E_j$ with respect to the partial order if $E_i$ is closer to the braid axis than $E_j$; that is, if $E_i$ lies between lower-numbered strands than $E_j$ does. Let $E_1$ denote the innermost elementary region, containing the braid axis. Label the other elementary regions $E_2,\ldots,E_m$ so that whenever $i<j$, $E_i$ is less than or not comparable to $E_j$ with respect to the partial order.

\begin{definition}[Coherent Regions]
\label{dfnregions}
A \emph{coherent region} in the $I$-resolution of $D$ is the union of a set of non-comparable elementary regions, along with all elementary regions less than these under the partial order described above.  The \emph{weight} $\emph{\weight}(R)$ of a coherent region $R$ is twice the number of 4-valent vertices plus the number of bivalent vertices in the closure of $R$.  The \emph{non-local relation associated to $R$} is $$t^{\emph{\weight}(R)}w_{\text{out}}-w_{\text{in}},$$ where $w_{\text{out}}$ (respectively $w_{\text{in}}$) is the product of the edges outside $R$, but incident to exactly one vertex of $\partial R$ and pointing out from (respectively into) $R$.
\end{definition}

There are five coherent regions in the singularized figure 8 example of Figure~\ref{fig:nonlocal}, with associated relations as follows. Notice that, for example, $E_1\cup E_2\cup E_4$ is not a coherent region because $E_3 < E_4$.
\smallskip
\begin{center}
\renewcommand\arraystretch{1.5}
\begin{tabular}{|l|l|}
\hline
coherent region & non-local relation\\\hline
$E_1$ & $t^6x_1x_7-x_4x_{10}$ \\\hline
$E_1\cup E_2$ & $t^8x_1x_9-x_4x_6$ \\\hline
$E_1\cup E_3$ & $t^8x_3x_7-x_0x_{10}$ \\\hline
$E_1\cup E_2\cup E_3$ & $t^{10}x_3x_9-x_0x_6$ \\\hline
$E_1\cup E_2\cup E_3\cup E_4$ & $t^{11}x_9-x_0$\\
\hline
\end{tabular}
\end{center}
\smallskip

Finally, we may think of non-local relations as arising from subsets of vertices in the $I$-resolution of $D$. 

\begin{definition}[Subsets]
\label{dfnsubsets}
Let $V$ be a subset of the vertices in the $I$-resolution of $D$. The \emph{weight} $\emph{\weight}(V)$ of $V$ is twice the number of 4-valent vertices plus the number of bivalent vertices in $V$. The \emph{non-local relation associated to $V$} is $$t^{\emph{\weight}(V)}w_{\text{out}}-w_{\text{in}},$$ where $w_{\text{out}}$ is the product of edges from $V$ to its complement and $w_{\text{in}}$ is the product of edges into $V$ from its complement.
\end{definition}

Any of these three definitions gives a generating set for $\cN_I(D)$. We will prove that the three definitions are equivalent in Proposition~\ref{nonlocaldfns}. First, we record some observations about the efficiency of the generating sets prescribed by the different definitions.

A priori, the generating set obtained from subsets is much larger than those obtained from cycles or coherent regions. However, it actually suffices to consider a smaller collection of subsets whose associated relations still generate the same ideal in $\cR[x_0,\ldots,x_n]/\cL$.  First, we may restrict to connected subsets of vertices, meaning those whose union with their incident edges is a connected graph. If a subset $V$ is disconnected as $V=V^\prime\coprod V^{\prime\prime}$, then the outgoing (respectively incoming) edges from $V$ are exactly the union of the outgoing (respectively incoming) edges from $V^\prime$ and $V^{\prime\prime}$. Therefore, the non-local relation associated to $V$ has the form $$t^{\weight(V^\prime)+\weight(V^{\prime\prime})}w_{\text{out}}^\prime w_{\text{out}}^{\prime\prime}-w_{\text{in}}^\prime w_{\text{in}}^{\prime\prime}.$$ However, this is already contained in the ideal generated by $$t^{\weight(V^\prime)}w_{\text{out}}^\prime-w_{\text{in}}^\prime\quad\text{and}\quad t^{\weight(V^{\prime\prime})}w_{\text{out}}^{\prime\prime}-w_{\text{in}}^{\prime\prime},$$ which are the non-local relations associated to $V^\prime$ and $V^{\prime\prime}$.

Second, we may ignore a subset $V$ if the union of $V$ with its incident edges is a graph with no oriented cycles. In Figure~\ref{fig:nonlocal}, the two vertices in layer $s_0$ along with the 4-valent vertex in layer $s_1$ form such a subset. The non-local relation associated to this subset is $t^5x_3x_7x_8-x_0x_1x_2$, but simple substitutions using the local relations associated to the three vertices in the subset show that this supposedly non-local relation is actually contained in $\cL$.

\begin{obs}
\label{subsetobs}
The ideal of non-local relations $\cN_I$ can be generated in $\cR[\underline{x}(D)]/\cL$ by the non-local relations associated to connected subsets that contain oriented cycles.
\end{obs}

We prove this statement inductively. The base case is a subset with a single vertex and no oriented cycles, which has associated non-local relation identical to its associated quadratic local relation. Any connected subset $V=\{v_1,\ldots,v_k\}$ with no oriented cycles can be constructed from a single vertex ($v_1$, renumbering as necessary) by adding vertices one at a time such that each $V_i=\{v_1,\ldots,v_i\}$ is a connected subset with no oriented cycles and $v_{i+1}$ is connected to $V_i$ only by edge(s) outgoing from $V_i$. Let $\cN_I^\prime$ be the ideal generated by non-local relations associated to connected subsets that contain oriented cycles. Assume that the non-local relation associated to $V_i$ is in the ideal sum $\cL+\cN_I^\prime$.

Suppose first that $v_{i+1}$ is a bivalent vertex. Let $x_{\text{out}}$ be the edge from $V_i$ to $v_{i+1}$ and $x_v$ is the edge pointing out from $v_{i+1}$. If $t^{\weight(V_i)} w_{\text{out}}x_{\text{out}}-w_{\text{in}}$ is the non-local relation associated to $V_i$, then the relation associated to $V_{i+1}$ is $t^{\weight(V_i) +1}w_{\text{out}}x_v-w_{\text{in}}$. Using the local relation $tx_v-x_{\text{out}}$ to replace $x_v$ recovers the non-local relation associated to $V_{i}$. So the non-local relation associated to $V_{i+1}$ is contained in the ideal sum of $\cL$ with the non-local relation associated to $V_i$, hence in $\cL+\cN_I^\prime$. 

Suppose instead that $v_{i+1}$ is a 4-valent vertex with edges $x_a$ and $x_b$ pointing out and edges $x_c$ and $x_{\text{out}}$ pointing in. Suppose that $x_{\text{out}}$ connects to a vertex in $V_i$ and that $x_c$ is not incident to any vertex in $V_i$. The local relation associated to $v_{i+1}$ is then $t^2x_ax_b-x_cx_{\text{out}}$, while the non-local relation associated to $V_i$ is of the form $t^{\weight (V_i)}w_{\text{out}}x_{\text{out}}-w_{\text{in}}$.  The non-local relation associated to $V_{i+1}$ is $$t^{\weight (V_i)+2}w_{\text{out}}x_ax_b-w_{\text{in}}x_c\equiv t^{\weight (V_i)}w_{\text{out}}x_cx_{\text{out}}-w_{\text{in}}x_c=x_c(t^{\weight (V_i)}w_{\text{out}}x_{\text{out}}-w_{\text{in}}).$$  Therefore, the non-local relation associated to $V_{i+1}$ is again contained in the ideal sum of $\cL$ with the non-local relation associated to $V_i$, hence in $\cL+\cN_I^\prime$.

Finally, it may be that $v_{i+1}$ is a 4-valent vertex connected to $V_i$ by two edges. Let $x_a$ and $x_b$ be the edges pointing out from $v_{i+1}$ and $x_{c}$ and $x_{d}$ the edges pointing into $v_{i+1}$. Then the quadratic relation associated to $v_{i+1}$ is $t^2x_ax_b-x_{c}x_{d}$. The non-local relation associated to $V_i$ is of the form $t^{\weight(V_i)}w_{\text{out}}x_{c}x_{d}-w_{\text{in}}$. The non-local relation associated to $V_{i+1}$ is \[t^{\weight(V_i)+2}w_{\text{out}}x_ax_b-w_{\text{in}}=t^{\weight(V_i)}w_{\text{out}}\left(t^2x_ax_b-x_cx_d\right)+t^{\weight(V_i)}w_{\text{out}}x_{c}x_{d}-w_{\text{in}}.\]

The second observation of this section concerns redundancy in the generating sets for $\cN_I$ defined by cycles and coherent regions arising from certain elementary regions that can be removed from a coherent region without producing an independent non-local relation.  For instance, in Figure~\ref{fig:nonlocal}, the coherent region $E_1\cup E_2\cup E_3\cup E_4$ specifies the non-local relation $t^{11}x_9-x_0$ as a generator for $\cN_I$. Then $x_6(t^{11}x_9-x_0)$ is also in $\cN_I$. It can be modified to $t^{10}x_3x_9-x_0x_6$ using the relation $tx_6-x_3$, which is the linear relation associated to the bivalent vertex in layer $s_1$.  We have obtained the non-local relation associated to $E_1\cup E_2\cup E_3$, showing that it is redundant once the non-local relation for $E_1\cup E_2\cup E_3\cup E_4$ is included in the generating set of $\cN_I$. More formally, we have the following observation.

\begin{obs}
\label{regionobs}
Suppose a coherent region $R$ has an adjacent elementary region $E$ and that $\partial E\setminus\partial R\cap\partial E$ is a path of edges through bivalent vertices only. Then the non-local relation associated to $R$ is contained in the ideal sum of $\cL$ with the non-local relation associated to $R\cup E$.
\end{obs}

Label the edges in the path in $\partial E\setminus\partial R\cap\partial E$ by $x_{\text{out}},x_1,\ldots,x_p,x_{\text{in}}$ consistent with the orientation of the overall diagram. The linear relations associated to each vertex in this path are $tx_1-x_{\text{out}}$, $tx_{\text{in}}-x_p$, and $tx_{i+1}-x_i$ for $1\leq i\leq p-1$. The non-local relation associated to $R$ has the form $$t^{\weight(R)} w_{\text{out}}x_{\text{out}}-w_{\text{in}}x_{\text{in}},$$ which can be rewritten using the linear relations above to give $$t^{\weight(R) +p+1}w_{\text{out}}x_{\text{in}}-w_{\text{in}}x_{\text{in}}=\left(t^{\weight(R) +p+1}w_{\text{out}}-w_{\text{in}}\right)x_{\text{in}},$$ which is a multiple of the non-local relation associated to $R\cup E$. Therefore, to form a minimal generating set for $\mathcal{N}_I$, we need only consider $R\cup E$.

As these observations begin to indicate, the definitions of non-local relations via cycles, coherent regions, and subsets are equivalent. In the example of Figure~\ref{fig:nonlocal}, the cycle shown in bold produces the same non-local relation as the coherent region $E_1\cup E_2$ or the subset of vertices contained in $E_1\cup E_2$. These correspondences between cycles, coherent regions, and subsets hold in general.

\begin{prop}
\label{nonlocaldfns}
Definitions~\ref{dfncycles}, \ref{dfnregions}, and \ref{dfnsubsets} produce the same ideal in $\cR[x_0,\ldots,x_n] / \cL$, where $\cL$ is the ideal generated by local relations associated to each vertex in the $I$-resolution of $D$.
\end{prop}
\begin{proof}
The equivalence between definitions~\ref{dfncycles} (cycles) and \ref{dfnregions} (coherent regions) is clear: the boundaries of coherent regions are exactly the cycles that avoid the basepoint and have orientations matching that of $D$. (Consider, for example, the boundary of $E_1\cup E_2\cup E_3$ compared to that of $E_1\cup E_2\cup E_4$ in Figure~\ref{fig:nonlocal}.) Weights and the edge products $w_{\text{out}}$ and $w_{\text{in}}$ are identical for a coherent region $R$ and the cycle $\partial R$, so the associated non-local relations are the same.

Let $\cN$ denote the ideal generated by non-local relations associated to cycles or coherent regions in $\cR[\underline{x}(D)]/\cL$. Let $\cN_S$ denote the ideal generated by non-local relations associated to subsets.  Suppose $R$ is a coherent region and $V_R$ the set of vertices in its closure.  Then $\weight (R)=\weight (V_R)=\weight (\partial R)$ and the words $w_{\text{out}}$ and $w_{\text{in}}$ defined with respect to $R$, $\partial R$, or $V_R$ are the same. Therefore, we have the inclusion $\cN\subset\cN_S$. 

For the opposite inclusion, consider a subset $V$.  We appeal first to Observation~\ref{subsetobs}, which allows us to assume that the union of $V$ and its incident edges forms a connected graph containing an oriented cycle $Z$. Assume that $Z$ is the outermost cycle contained in $V$, and let $R_Z$ be the coherent region it bounds. If $V$ contains all of the vertices in the closure of $R_Z$, then the arguments about connectedness and subsets just after Observation~\ref{subsetobs} allow us to remove all vertices from $V$ that are not contained in the closure of $R_Z$, thereby showing that the non-local relation associated to $V$ can be constructed from the non-local relation associated to $R_Z$.

Suppose now that $V$ does not contain all of the vertices in $R_Z$. Then the complement of $V$ is disconnected, with one component inside $Z$ and one component outside. Denote these components $V^\prime$ and $V^{\prime\prime}$, respectively. Then $V\cup V^\prime$ contains $Z$ and all of the vertices in the closure of $R_Z$, so the argument above shows that its associated non-local relation is contained in $\cN$. The subset $V^\prime$ may not contain any oriented cycles or it may contain an oriented cycle $Z^\prime$ and all vertices in the closure of $R_{Z^\prime}$. Therefore, its associated non-local relation is contained in either $\cL$ or $\cN$.

Finally, we show that the non-local relation associated to $V$ is in the ideal generated by the non-local relations associated to $V^\prime$ and $V\cup V^\prime$. The words $w_{\text{out}}$ and $w_{\text{in}}$ defined with respect to $V$ are products $w_{\text{out}}=w_{\text{in}}^\prime w_{\text{in}}^{\prime\prime}$ and $w_{\text{in}}=w_{\text{out}}^\prime w_{\text{out}}^{\prime\prime}$ of edges into and out from $V^\prime$ and $V^{\prime\prime}$. 
\begin{align*}
&t^{\weight (V)}w_{\text{in}}^\prime w_{\text{in}}^{\prime\prime}-w_{\text{out}}^\prime w_{\text{out}}^{\prime\prime} & \text{non-local relation from $V$}\\
\equiv& \,\,t^{\weight (V)+\weight (V^\prime)}w_{\text{out}}^\prime w_{\text{in}}^{\prime\prime}-w_{\text{out}}^\prime w_{\text{out}}^{\prime\prime} &\text{by substituting non-local relation from $V^\prime$}\\
=&\,\,(t^{\weight (V)+\weight (V^{\prime})}w_{\text{in}}^{\prime\prime}- w_{\text{out}}^{\prime\prime})w_{\text{out}}^\prime &\text{a multiple of the non-local relation $V\cup V^\prime$}
\end{align*}
Since the non-local relation associated to $V$ can be constructed from those associated to $V^\prime$ and $V\cup V^\prime$, it is contained in $\cN$. Therefore, any non-local relation associated to a subset can be generated from non-local relations associated to coherent regions, meaning that $\cN_S\subset\cN$.
\end{proof}

Three further observations related to the non-local relations are worth recording for later use. 

\begin{obs}
\label{singsmoothobs}
Let $D^\times$ and $D^\Arrowvert$ be resolutions of a diagram in which a certain crossing is singularized and smoothed, respectively, and that are otherwise identical. Then $\cN(D^\times)\subset\cN(D^\Arrowvert)$. 
\end{obs}

Let $v$ be the 4-valent vertex in $D^\times$ corresponding to the given crossing and $v_1$, $v_2$ the two bivalent vertices in $D^\Arrowvert$ corresponding to the same crossing. Let $V$ be a subset of the vertices in $D^\times$ that does not contain $v$. Then the non-local relation associated to $V$ is the same in $D^\times$ and in $D^\Arrowvert$. The non-local relation associated to $V\cup\left\{v\right\}$ in $D^\times$ is the same as the non-local relation associated to $V\cup\left\{v_1,v_2\right\}$ in $D^\Arrowvert$. Therefore, $\cN(D^\times)\subset\cN(D^\Arrowvert)$, as claimed. The opposite inclusion is false in general because $V\cup\left\{v_1\right\}$ and $V\cup\left\{v_2\right\}$ may have associated relations in $\cN(D^\Arrowvert)$ that are not contained in $\cN(D^\times)$.

\begin{obs}
\label{outermostcycleobs}
Let $\emph{{\weight}}(D)$ mean the weight of the set of all vertices in any resolution of $D$. This weight does not depend on the resolution of $D$ under consideration. The relation $t^{\emph{\weight} (D)}x_n-x_0$, where $x_n$ is the edge entering the basepoint and $x_0$ is the edge leaving it, holds in $\cA_I(D)$ for any $I$ and any $D$. It is associated to the subset containing all vertices or the outermost cycle in the diagram that does not pass through the basepoint.
\end{obs}

In a sense, then, the $\ast$ behaves like a bivalent vertex with weight $-\weight(D)$, balancing out the weight in the rest of the diagram.

\begin{obs}
\label{disconnectedobs}
If $I$ is a disconnected resolution of $D$, and we choose to work over a completed ground ring, then the algebra associated to the $I$-resolution of $D$ will vanish. In a disconnected resolution, there are cycles that do not contain the basepoint and have no ingoing or outgoing edges. In this situation, we interpret the products $w_\text{out}$ and $w_\text{in}$ to be $1$, which makes the associated relation $t^k-1$ for some $k$. In $\widehat{\cR}$ or $\widehat{\cR[\underline{x}]}$, $t^k-1$ is a unit. Therefore, including $t^k-1$ in our ideal of relations makes $\cA_I(D)\otimes_{\cR}\widehat{\cR}$ or $\cA_I(D)\otimes_{\cR[\underline{x}]}\widehat{\cR[\underline{x}]}$ vanish.
\end{obs}

The fact that $\cA_I(D)=0$ for any disconnected resolution of $D$ is a significant distinction between this construction and the HOMFLY-PT homology of \cite{kr2}, but is appropriate in light of the fact that the Alexander polynomial vanishes on split links.

\section{Removing bivalent vertices}
\label{sec:removemark}

\begin{figure}[htbp]
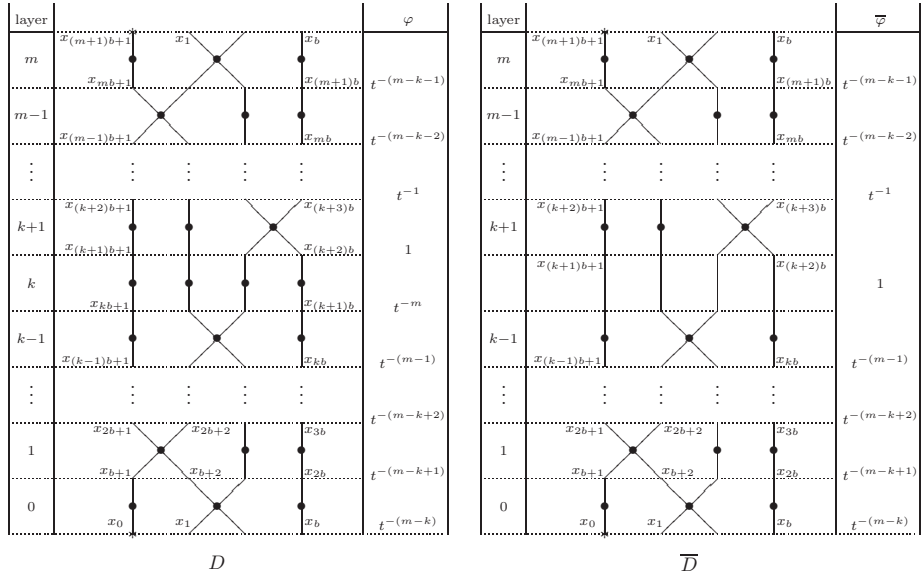

	\begin{center}
			\subfloat{\scalebox{.74}{\input{removemarklayer.tex}}}\quad
			\subfloat{\scalebox{.74}{\input{removemarkdone.tex}}}
	\caption{Diagrams for the proof of Proposition~\ref{prop:removemark}. The maps $\varphi$ on $\cA_I(D)$ and $\overline{\varphi}$ on $\cA_{\overline{I}}(\overline{D})$ are defined to be multiplication by the factor shown in the rightmost column of each diagram.}
	\label{removemarkfig}
	\end{center}
\end{figure}

This section is devoted to a technical result allowing us to remove a horizontal layer of a diagram with a bivalent vertex on each strand and no 4-valent vertices. Such a layer is obtained each time a crossing is resolved. Suppose the $I$-resolution of $D$ is a diagram with $m+1$ layers, and that layer $k$ contains only bivalent vertices.  Let $\overline{D}$ denote the diagram obtained by removing layer $k$. The proposition below shows that adding or removing layer $k$ corresponds to twisting the action of the ground ring via a non-trivial endomorphism. We may either describe $\cA_I(D)$ and $\cA_{\overline{I}}(\overline{D})$ each as twistings of a common $\cR[\underline{x}(\overline{D})]$-module, or we may describe them as twistings of each other after enlarging the ground ring to include appropriate roots of $t$. Applying any of these twistings to every summand of a chain complex $C(D)$ would not change the homology of the complex because $\cR$ is flat when considered as an $\cR$-module via any of the relevant endomorphisms. We refer to the notation in Figure~\ref{removemarkfig} throughout.

\begin{prop}
\label{prop:removemark}
Let $D$ and $\overline{D}$ be defined as above.  Let $\overline{I}$ denote the index $I$ with its $k^{th}$ component deleted. Then there is an $\cR[\underline{x}(\overline{D})]$-module $\cA_{I,\overline{I}}(D,\overline{D})$ and $\cR[\underline{x}(\overline{D})]$-module isomorphisms 
\begin{eqnarray*}
\cA_I(D)&\isom&\cA_{I,\overline{I}}(D,\overline{D})\otimes_{(\cR,\psi_{m+1})}\cR\\
\cA_{\overline{I}}(\overline{D})&\isom&\cA_{I,\overline{I}}(D,\overline{D})\otimes_{(\cR,\psi_{m})}\cR,
\end{eqnarray*}
where $\psi_i$ is the endomorphism of $\cR$ taking $1$ to $1$ and $t$ to $t^{i}$.

Let $\cA_I^{(i)}(D)=\cA_I(D)\otimes_\cR\cR[t^{1/i}]$ and $\cA_{\overline{I}}^{(i)}(\overline{D})=\cA_{\overline{I}}(\overline{D})\otimes_\cR\cR[t^{1/i}]$. There is an $\cR[t^{1/m}][\underline{x}(\overline{D})]$-module isomorphism
\begin{equation*}
\cA^{(m)}_I(D)\isom\cA^{(m)}_{\overline{I}}(\overline{D})\otimes_{(\cR[t^{1/m}],\psi_{m+1}\circ\psi_m^{-1})}\cR[t^{1/m}]\\
\end{equation*}
and an $\cR[t^{1/(m+1)}][\underline{x}(\overline{D})]$-module isomorphism
\begin{equation*}
\cA^{(m+1)}_{\overline{I}}(\overline{D})\isom\cA^{(m+1)}_{I}(D)\otimes_{(\cR[t^{1/(m+1)}],\psi_{m}\circ\psi_{m+1}^{-1})}\cR[t^{1/(m+1)}].
\end{equation*}
\end{prop}

\begin{proof}
We first define automorphisms $\varphi$ of $\cA_I(D)$ and $\overline{\varphi}$ of $\cA_I(\overline{D})$ that transform our original presentations of these algebras into presentations in which $t$ appears very rarely.

Define $\varphi$ to be multiplication by $t^{-(j-1)}$ on edges $x_{(k+j)b+i}$ for $0\leq j\leq m$ and $1\leq i\leq b$ (treating the $k+j$ portion of the subscript modulo $m+1$), and multiplication by $t^{-(m-k)}$ on edge $x_{(m+1)b+1}=x_n$. That is, $\varphi$ is the identity on the edges connecting layer $k$ to layer $k+1$ (edges $x_{(k+1)b+1},\ldots,x_{(k+2)b}$), multiplication by $t^{-1}$ on the edges connecting layer $k+1$ to layer $k+2$ (edges $x_{(k+2)b+1},\ldots,x_{(k+3)b}$), multiplication by $t^{-2}$ on the edges connecting layer $k+2$ to layer $k+3$, and so on, until it is multiplication by $t^{-m}$ on the edges connecting layer $k-1$ to layer $k$ (edges $x_{kb+1},\ldots,x_{(k+1)b}$).

We may continue to use $x_0,\ldots,x_n$ as generators of $\varphi(\cA_I(D))$, but must examine carefully the effect of $\varphi$ on the generating sets of $\cL$ and $\cN_I(D)$. Consider first the generators of $\cL_{k+j}$ for any $j\neq 0$. These have one of the following forms, where $1\leq i\leq b$.
\begin{align}
tx_{(k+j+1)b+i}+tx_{(k+j+1)b+i+1}&-x_{(k+j)b+i}-x_{(k+j)b+i+1}\\
t^2x_{(k+j+1)b+i}x_{(k+j+1)b+i+1}&-x_{(k+j)b+i}x_{(k+j)b+i+1}\\
tx_{(k+j+1)b+i}&-x_{(k+j)b+i}
\end{align}
After applying $\varphi$, they become 
\begin{align}
\label{represent1}
t^{-j+1}x_{(k+j+1)b+i}+t^{-j+1}x_{(k+j+1)b+i+1}&-t^{-(j-1)}x_{(k+j)b+i}-t^{-(j-1)}x_{(k+j)b+i+1}\nonumber\\
\equiv x_{(k+j+1)b+i}+x_{(k+j+1)b+i+1}&-x_{(k+j)b+i}-x_{(k+j)b+i+1}
\end{align}
\begin{align}
t^{-2j+2}x_{(k+j+1)b+i}x_{(k+j+1)b+i+1}&-t^{-2(j-1)}x_{(k+j)b+i}x_{(k+j)b+i+1}\nonumber\\
\equiv x_{(k+j+1)b+i}x_{(k+j+1)b+i+1}&-x_{(k+j)b+i}x_{(k+j)b+i+1}
\end{align}
\begin{align}
\label{represent3}
t^{-j+1}x_{(k+j+1)b+i}&-t^{-(j-1)}x_{(k+j)b+i}\nonumber\\
\equiv x_{(k+j+1)b+i}&-x_{(k+j)b+i}.
\end{align}

The price of eliminating powers of $t$ from most local relations is that $t$ appears with higher powers in relations that do involve layer $k$. Since layer $k$ has only bivalent vertices, its associated relations are all of the form $tx_{(k+1)b+i}-x_{kb+i}$. Applying $\varphi$, we obtain $$tx_{(k+1)b+i}-t^{-m}x_{kb+i} \equiv t^{m+1}x_{(k+1)b+i}-x_{kb+i}.$$ 

Non-local relations are similarly affected. Consider the generating set for $\cN_I(D)$ given by coherent regions. We will show that $\varphi$ applied to any relation in this generating set produces a relation of the form $t^{p(m+1)}w_{\text{out}}-w_{\text{in}}$ for some integer $p$. Begin with the innermost elementary region $E_1$. Suppose it has $v$ 4-valent vertices along its boundary in layers $k+j_1,\ldots,k+j_v$. Then $\weight (E_1)=m+1+v$.  Each 4-valent vertex contributes one edge to the product $w_{\text{out}}$ and an edge one layer lower to $w_{\text{in}}$. If $j_i\neq 0$ for $1\leq i\leq v$, then
\begin{align*}
&\varphi\left(w_{\text{out}}\right)=t^{-j_1-\cdots-j_v}w_{\text{out}}&\text{and}\\
&\varphi\left(w_{\text{in}}\right)=t^{-(j_1-1)-\cdots-(j_v-1)}w_{\text{in}}=t^{-j_1-\cdots-j_v+v}w_{\text{in}},&\text{so}\\
&\varphi\left(t^{m+1+v}w_{\text{out}}-w_{\text{in}}\right)\equiv t^{m+1}w_{\text{out}}-w_{\text{in}}.
\end{align*}
Suppose instead (without loss of generality) that $j_1=0$. Then $\varphi$ is the identity when applied to the outgoing edge of the vertex in layer $k+j_1$, but multiplication by $t^{-m}$ on the incoming edge. Therefore, 
\begin{align*}
&\varphi\left(w_{\text{out}}\right)=t^{-j_2-\cdots-j_v}w_{\text{out}}&\text{and}\\
&\varphi\left(w_{\text{in}}\right)=t^{-m-(j_2-1)-\cdots-(j_v-1)}w_{\text{in}}=t^{-m-j_2-\cdots-j_v+v-1}w_{\text{in}},&\text{so}\\
&\varphi\left(t^{m+1+v}w_{\text{out}}-w_{\text{in}}\right)\equiv t^{2(m+1)}w_{\text{out}}-w_{\text{in}}.
\end{align*}
So $\varphi$ has the claimed effect on the non-local relation associated to the innermost coherent region.

Next consider an elementary region $E\neq E_1$ with bottommost vertex in layer $k+j$ and topmost vertex in layer $k+j+s$. Suppose $\partial E$ meets $v^\prime$ additional 4-valent vertices in layers $k+j_1,\ldots,k+j_{v^\prime}$.  Assume for now that $E$ does not meet layer $k$. Then $\weight(E)=2(s+1)+v^\prime$. Let $t^{2(s+1)+v^\prime}e_{\text{out}}-e_{\text{in}}$ denote the non-local relation associated to $E$. The topmost vertex of $E$ contributes two outgoing edges to $e_{\text{out}}$ and the bottommost vertex contributes two incoming edges to $e_{\text{in}}$.  The other $v^\prime$ 4-valent vertices contribute one edge each to $e_{\text{out}}$ and $e_{\text{in}}$. Therefore, 
\begin{align*}
&\varphi\left(e_{\text{out}}\right)=t^{-2(j+s)-j_1-\cdots-j_{v^\prime}}e_{\text{out}}=t^{-2j-j_1-\cdots-j_{v^\prime}-2s}e_{\text{out}}&\text{and}\\
&\varphi\left(e_{\text{in}}\right)=t^{-2(j-1)-(j_1-1)-\cdots-(j_{v^\prime}-1)}e_{\text{in}}=t^{-2j-j_1-\cdots-j_{v^\prime}+v^\prime+2}e_{\text{in}},&\text{so}\\
&\varphi\left(t^{2(s+1)+v^\prime}e_{\text{out}}-e_{\text{in}}\right)\equiv e_{\text{out}}-e_{\text{in}}.
\end{align*}
\noi If $E$ does meet layer $k$, a then modification of the calculation above (similar to that used for $E_1$) verifies the claim that $\varphi\left(t^{2(s+1)+v^\prime}e_{\text{out}}-e_{\text{in}}\right)$ has the form $t^{p(m+1)}e_{\text{out}}-e_{\text{in}}$ for some integer $p$.

Finally, consider a coherent region $R^\prime$ that is not elementary. We can write $R^\prime$ as $R\cup E$, where $R$ is a coherent region and $E$ is an elementary region. Suppose the non-local relations associated to $R$ and $E$ are
$t^{\weight(R)}w_{\text{out}}-w_{\text{in}}$ and $t^{\weight(E)}e_{\text{out}}-e_{\text{in}},$ respectively. Let $y$ be the product of edges that connect vertices in $R$ to vertices in $E$. The non-local relation associated to $R^\prime$ can be obtained by combining the non-local relations associated to $R$ and $E$, then factoring out $y$ as follows.
$$t^{\weight(R)+\weight(E)}w_{\text{out}}e_{\text{out}}-w_{\text{in}}e_{\text{in}}=y\left(t^{\weight(R)+\weight(E)}w^\prime_{\text{out}}e^\prime_{\text{out}}-w^\prime_{\text{in}}e^\prime_{\text{in}}\right)$$
The non-local relation associated to $R^\prime$ is $t^{\weight(R)+\weight(E)}w^\prime_{\text{out}}e^\prime_{\text{out}}-w^\prime_{\text{in}}e^\prime_{\text{in}}$. We will assume inductively that $\varphi$ applied to the non-local relations for $R$ and $E$ produces $t^{p(m+1)}w_{\text{out}}-w_{\text{in}}$ and $t^{q(m+1)}e_{\text{out}}-e_{\text{in}}$, respectively for some integers $p$ and $q$.  Then 
\begin{align*}
&\varphi\!\left(t^{\weight(R)+\weight(E)}w_{\text{out}}e_{\text{out}}-w_{\text{in}}e_{\text{in}}\right)\\
\equiv\,&t^{(p+q)(m+1)}w_{\text{out}}e_{\text{out}}-w_{\text{in}}e_{\text{in}}\\
=\,&y\!\left(t^{(p+q)(m+1)}w^\prime_{\text{out}}e^\prime_{\text{out}}-w^\prime_{\text{in}}e^\prime_{\text{in}}\right)
\end{align*}
and on the other hand
\begin{align*}
&\varphi\!\left(t^{\weight(R)+\weight(E)}w_{\text{out}}e_{\text{out}}-w_{\text{in}}e_{\text{in}}\right)\\
\equiv\,& \varphi\!\left(y\right)\varphi\!\left(t^{\weight(R)+\weight(E)}w^\prime_{\text{out}}e^\prime_{\text{out}}-w^\prime_{\text{in}}e^\prime_{\text{in}}\right)\\
\equiv\,& y\,\varphi\!\left(t^{\weight(R)+\weight(E)}w^\prime_{\text{out}}e^\prime_{\text{out}}-w^\prime_{\text{in}}e^\prime_{\text{in}}\right).
\end{align*}
We have verified that applying $\varphi$ to the non-local relation associated to $R^\prime$ produces a relation in which the power of $t$ is an integer multiple of $m+1$.

So far, we have relations of the following forms in our presentation of $\varphi(\cA_I(D))$.
\begin{align}
x_{(k+j+1)b+i}+x_{(k+j+1)b+i+1}&-x_{(k+j)b+i}-x_{(k+j)b+i+1}\nonumber\\
x_{(k+j+1)b+i}x_{(k+j+1)b+i+1}&-x_{(k+j)b+i}x_{(k+j)b+i+1}\nonumber\\
x_{(k+j+1)b+i}&-x_{(k+j)b+i}\nonumber\\
t^{m+1}x_{(k+1)b+i}&-x_{kb+1}\label{useme}\\
t^{p(m+1)}w_{\text{out}}&-w_{\text{in}}\nonumber
\end{align}

\noindent It will be convenient to make one final modification: use the relations in (\ref{useme}) to eliminate the variables for edges connecting layer $k-1$ to layer $k$. The result is a presentation in which $t$ appears only in the following types of relations.

\begin{align}
\label{match1}
t^{m+1} x_{(k+1)b+i}+t^{m+1} x_{(k+1)b+i+1}&-x_{(k-1)b+i}-x_{(k-1)b+i+1}\\
\label{match2}
t^{2(m+1)}x_{(k+1)b+i}x_{(k+1)b+i+1}&-x_{(k-1)b+i}x_{(k-1)b+i+1}\\
\label{match3}
t^{m+1}x_{(k+1)b+i}&-x_{(k-1)b+i}\\
\label{match4}
t^{p(m+1)} w_{\text{out}}&-w_{\text{in}}
\end{align}

The second map, $\overline{\varphi}$, allows us to present $\cA_{\overline{I}}(\overline{D})$ in a similar way, with powers of $t$ appearing only in certain relations, and only as $t^{pm}$ for various integers $p$. Define $\overline{\varphi}$ in exactly the same way as  $\varphi$ on edges $x_{(k+j)b+i}$ for $1\leq j\leq m$ and $0\leq i \leq b$ and for edge $x_{(m+1)b+1}$. Diagram $\overline{D}$ has no $k^{th}$ layer, so $\overline{\varphi}$ is the identity on the edges connecting layer $k-1$ to layer $k+1$, multiplication by $t^{-1}$ on the edges connecting layer $k+1$ to layer $k+2$, multiplication by $t^{-2}$ on the edges connecting layer $k+2$ to $k+3$, and so on, until it is multiplication by $t^{-(m-1)}$ on the edges connecting layer $k-2$ to layer $k-1$.

Again, for most relations, $\overline{\varphi}$ eliminates all powers of $t$. Similar calculations to those above show that $\overline{\varphi}$ removes $t$ from the generating set for $\cL_{k+j}$ for $j\neq 0$, leaving relations identical to those in (\ref{represent1}) to (\ref{represent3}) above.

All powers of $t$ end up in generators of $\cL_{k-1}$ and $\cN_I$, but this time with multiples of $m$ instead of $m+1$. The relations that involve $t$ have one of the following forms.
\begin{align}
\label{match5}
t^{m} x_{(k+1)b+i}+t^{m} x_{(k+1)b+i+1}&-x_{(k-1)b+i}-x_{(k-1)b+i+1}\\
\label{match6}
t^{2m}x_{(k+1)b+i}x_{(k+1)b+i+1}&-x_{(k-1)b+i}x_{(k-1)b+i+1}\\
\label{match7}
t^{m}x_{(k+1)b+i}&-x_{(k-1)b+i}\\
\label{match8}
t^{pm} w_{\text{out}}&-w_{\text{in}}
\end{align}

The calculation proceeds as follows, where $p$ is the number of edges connecting layer $k-1$ used by a subset.
\begin{align}
\overline{\varphi}(tx_{(k+1)b+i}+tx_{(k+1)b+i+1}&-x_{(k-1)b+i}-x_{(k-1)b+i+1})\nonumber\\
=tx_{(k+1)b+i}+tx_{(k+1)b+i+1}&-t^{-(m-1)}x_{(k-1)b+i}-t^{-(m-1)}x_{(k-1)b+i+1}\nonumber\\
\equiv t^m x_{(k+1)b+i}+t^m x_{(k+1)b+i+1}&-x_{(k-1)b+i}-x_{(k-1)b+i+1}
\label{match5}
\end{align}
\begin{align}
\overline{\varphi}(t^2x_{(k+1)b+i}x_{(k+1)b+i+1}&-x_{(k-1)b+i}x_{(k-1)b+i+1})\nonumber\\
=t^2x_{(k+1)b+i}x_{(k+1)b+i+1}&-t^{-2(m-1)}x_{(k-1)b+i}x_{(k-1)b+i+1}\nonumber\\
\equiv t^{2m}x_{(k+1)b+i}x_{(k+1)b+i+1}&-x_{(k-1)b+i}x_{(k-1)b+i+1}
\label{match6}
\end{align}
\begin{align}
\overline{\varphi}(tx_{(k+1)b+i}&-x_{(k-1)b+i})\nonumber\\
=tx_{(k+1)b+i}&-t^{-(m-1)}x_{(k-1)b+i}\nonumber\\
\equiv t^m x_{(k+1)b+i}&-x_{(k-1)b+i}
\label{match7}
\end{align}
\begin{align}
\overline{\varphi}(t^{2s}w_{\text{out}}&-w_{\text{in}})\nonumber\\
\equiv t^{m p} w_{\text{out}}&-w_{\text{in}}
\label{match8}
\end{align}

We now have presentations of $\cA_I(D)$ and $\cA_{\overline{I}}(\overline{D})$, both over the smaller edge ring $\cR[\underline{x}(\overline{D})]$, that differ only by whether $t$ appears with a power of $m+1$ or with $m$. Here are two ways to describe the relationship between $\cA_I(D)$ and $\cA_{\overline{I}}(\overline{D})$. First, let $\cA_{I,\overline{I}}(D,\overline{D})$ be the $\cR[\underline{x}(\overline{D})]$-algebra obtained by setting $m=1$ in the presentation of $\overline{\varphi}(\cA_{\overline{I}}(\overline{D}))$ above (or, equivalently, setting $m+1=1$ in the presentation of $\varphi(\cA_I(D))$ above). Let $\psi_i$ be the endomorphism of $\cR$ that sends 1 to 1 and $t$ to $t^i$. Then $\cA_I(D)$ and $\cA_{\overline{I}}(\overline{D})$ are each obtained from $\cA_{I,\overline{I}}(D,\overline{D})$ by twisting by an endomorphism of $\cR$:
\begin{eqnarray*}
\cA_I(D)&\isom&\cA_{I,\overline{I}}(D,\overline{D})\otimes_{(\cR,\psi_{m+1})}\cR\\
\cA_{\overline{I}}(\overline{D})&\isom&\cA_{I,\overline{I}}(D,\overline{D})\otimes_{(\cR,\psi_{m})}\cR.
\end{eqnarray*}

Alternatively, we may enlarge the ground ring sufficiently to make $\psi_m$ or $\psi_{m+1}$ (or both) invertible and describe $\cA_I(D)$ and $\cA_{\overline{I}}(\overline{D})$ over the enlarged ground rings as twistings of each other rather than twistings of a third algebra. It follows directly from the preceding paragraph that:
\begin{eqnarray*}
\cA^{(m)}_I(D)&\isom&\cA^{(m)}_{\overline{I}}(\overline{D})\otimes_{(\cR[t^{1/m}],\psi_{m+1}\circ\psi_m^{-1})}\cR[t^{1/m}]\\
\cA^{(m+1)}_{\overline{I}}(\overline{D})&\isom&\cA^{(m+1)}_{I}(D)\otimes_{(\cR[t^{1/(m+1)}],\psi_{m}\circ\psi_{m+1}^{-1})}\cR[t^{1/(m+1)}].
\end{eqnarray*}
\end{proof}

\section{Braid-like Reidemeister Move II}
\label{sec:reid2}
\begin{figure}[htbp]
\begin{center}
\input{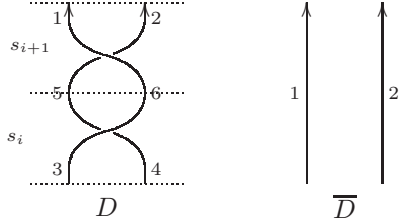}
\caption{Projection $D$ layers $s_i$ and $s_{i+1}$ and the corresponding portion of $\overline{D}$, which has no vertices. Technically, $\overline{D}$ does not have layers corresponding to $s_i$ and $s_{i+1}$; it is identical to $D$ in all other layers.  Assume that the braid axis is to the right of each diagram.}
\label{reid2labels}
\end{center}
\end{figure}

Suppose $D$ and $\overline{D}$ are two knot projections that differ by a Reidemeister II move with  labels as in Figure~$\ref{reid2labels}$.  The edge rings of $D$ and $\overline{D}$ are related by $\cR[\underline{x}(D)]=\cR[\underline{x}(\overline{D})][x_3,x_4,x_5,x_6]$. We will show that $C(D)$ and $C(\overline{D})$ are chain homotopy equivalent as complexes of $\cR[\underline{x}(\overline{D})]$-algebras, but will work over the larger edge ring $\cR[\underline{x}(D)]$ for as long as possible. Throughout this section, we will abbreviate indices of resolutions to two entries, showing only the states of the crossings in layers $s_i$ and $s_{i+1}$. 

\begin{figure}[htbp]
\begin{center}
\input{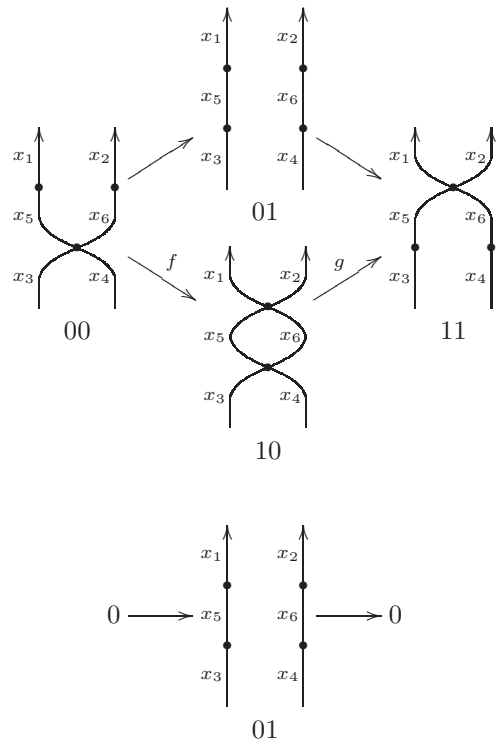}
\caption{The top chain complex is a portion of $C(D)$.  Lemma~\ref{lemma:reid2} shows that it is chain homotopy equivalent to the bottom chain complex. Assume that the braid axis is to the right of each diagram.}
\label{reid2resolutions}
\end{center}
\end{figure}

There are two oriented Reidemeister II moves, depending on which crossing in $D$ is positive and which is negative, but the arguments are very similar in the two cases. The relevant portion of $C(D)$ is shown in Figure~\ref{reid2resolutions}. The two variants of the Reidemeister II move exchange $\mathcal{A}_{00}(D)$ with $\mathcal{A}_{11}(D)$ and $\mathcal{A}_{01}(D)$ with $\mathcal{A}_{10}(D)$.

The key step in proving that the chain homotopy type of $C(D)$ is unchanged by a Reidemeister II move is to show the equivalence of the two complexes in Figure~$\ref{reid2resolutions}$. It suffices to prove the statement in Lemma~\ref{lemma:reid2}, which is a categorification of the MOY relation in the top line of Figure~\ref{fig:braidmoy}. Applying the direct sum splitting and an appropriate isomorphism of complexes, $\mathcal{A}_{00}(D)\stackrel{f}{\longrightarrow}\mathcal{A}_{10}(D)\stackrel{g}{\longrightarrow}\mathcal{A}_{11}(D)$ becomes an acyclic subcomplex.  Removing that subcomplex leaves the bottom complex of Figure~$\ref{reid2resolutions}$.  Removing the bivalent vertices in layers $s_i$ and $s_{i+1}$ (applying Proposition~\ref{prop:removemark} and reverting to the edge ring $\cR[\underline{x}(\overline{D})]$) leaves the corresponding portion of $C(\overline{D})$. 

\begin{lemma}
\label{lemma:reid2}
As $\cR[\underline{x}(\overline{D})][x_3,x_4]$-modules, $\mathcal{A}_{10}(D)\isom\mathcal{A}_{00}(D)\oplus\mathcal{A}_{11}(D)$, $f$ is an isomorphism onto the first summand, and $g$ is an isomorphism when restricted to the second summand.
\end{lemma}

\begin{proof}
The following matrix is a generating set for $\cL_i+\cL_{i+1}$ in the 10-resolution of $D$.
$$
\left(\begin{array}{c}
t(x_1+x_2)-(x_5+x_6)\\
(tx_1-x_6)(tx_2-x_6)\\
t(x_5+x_6)-(x_3+x_4)\\
(tx_6-x_4)(x_3-tx_6)
\end{array}\right)
$$
Use row I to eliminate $x_5$, then rewrite to limit the appearance of $x_6$ to a single row.
\begin{equation*}
{\left(\begin{array}{c}
(tx_1-x_6)(tx_2-x_6)\\
t^2(x_1+x_2)-(x_3+x_4)\\
(tx_6-x_4)(x_3-tx_6)
\end{array}\right)}
\xrightarrow{\text{III}+t^2\text{I}+tx_6\text{II}}
{\left(\begin{array}{c}
(tx_1-x_6)(tx_2-x_6)\\
t^2(x_1+x_2)-(x_3+x_4)\\
t^4x_1x_2-x_3x_4
\end{array}\right)}
\end{equation*}
Let $\overline{\cL}$ denote the ideal generated by the last two rows of the matrix above and $\cL$ denote the ideal generated by local relations in layers other than $i$ and $i+1$. Note that $x_5$ and $x_6$ do not appear in the generating set for $\cL$. By Observation~\ref{regionobs}, they need not appear in a generating set for $\cN_{10}$ either. Therefore, these ideals survive the manipulations above unchanged.  Define $$\cS=\frac{\cR[x_0,\ldots,x_4,x_7,\ldots,x_n]}{\overline{\cL}+\cL+\cN_{10}}.$$ We have simplified the presentation of $\cA_{10}$ so that $x_6$ appears only in one relation, which is quadratic in $x_6$. Using that relation, we may split $\cA_{10}$ as follows. $$\cA_{10}(D)\isom \frac{\cS[x_6]}{(tx_1-x_6)(tx_2-x_6)}\isom\cS(1)\oplus\cS(tx_1-x_6)$$ It remains to show that these two summands correspond to $\cA_{11}(D)$ and $\cA_{00}(D)$.

In the 11-resolution, the linear relations $tx_5-x_3$ and $tx_6-x_4$ may be used to replace $x_5$ and $x_6$ throughout the presentation. The resulting local relations in layers $i$ and $i+1$ exactly match those in $\overline{\cL}$. The definition by coherent regions and Observation~\ref{regionobs} give matching generating sets for $\cN_{10}$ and $\cN_{11}$. Therefore, $\cA_{11}(D)$ has a presentation identical to that of $\cS$ given above. Since $g$ is defined to be the quotient map, it is an isomorphism when restricted to the first summand of $\cA_{10}(D)$ above.

Similarly, in the 00-resolution, the linear relations $tx_1-x_5$ and $tx_2-x_6$ can be used to replace $x_5$ and $x_6$ throughout the presentation of $\cA_{00}(D)$. For local relations in layers $i$ and $i+1$, the resulting ideal is exactly $\overline{\cL}$. For non-local relations, the definition by coherent regions along with Observation~\ref{regionobs} again gives the same generating set for $\cN_{00}$ as for $\cN_{10}$. Therefore, $\cA_{00}(D)$ has a presentation identical to $\cS$. Since $f$ is defined to be multiplication by $tx_1-x_6$, it is an isomorphism onto the second summand of $\cA_{10}(D)$ above.
\end{proof}

\section{Braid-like Reidemeister Move III}
\label{sec:reid3}

\begin{figure}[thbp]
\begin{center}
\input{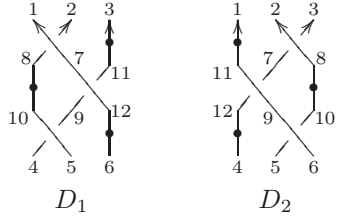}
\caption{Diagrams $D_1$ and $D_2$ in layers $s_1$, $s_2$, and $s_3$.}
\label{reid3adddots}
\end{center}
\end{figure}

In this section, we will consider projections $D_1$ and $D_2$ that differ by a Reidemeister III move with all negative crossings and labeling as in Figure~\ref{reid3adddots}.  Invariance under the other braid-like versions of Reidmeister III follows because all such moves are compositions of the negative Reidemeister III move and Reidemeister II moves.

\begin{figure}[phtb]
\begin{center}
\input{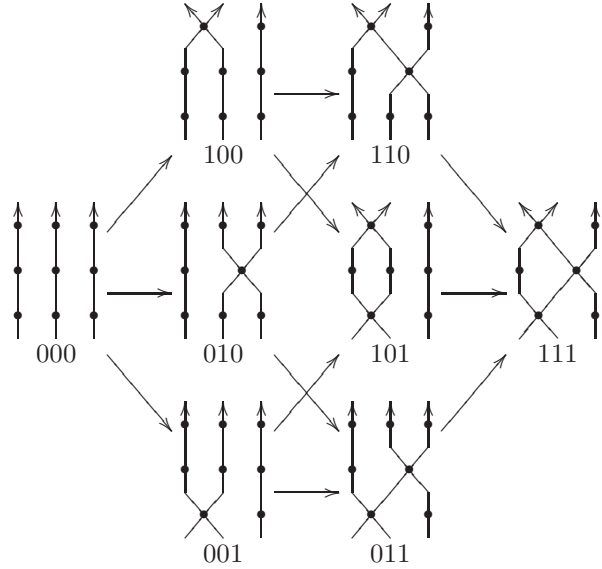}
\caption{Portion of the cube of resolutions for $D_1$ in layers $s_1$, $s_2$, and $s_3$.}
\label{reid3complexesneg}
\end{center}
\end{figure}

\begin{figure}[phtb]
\begin{center}
\input{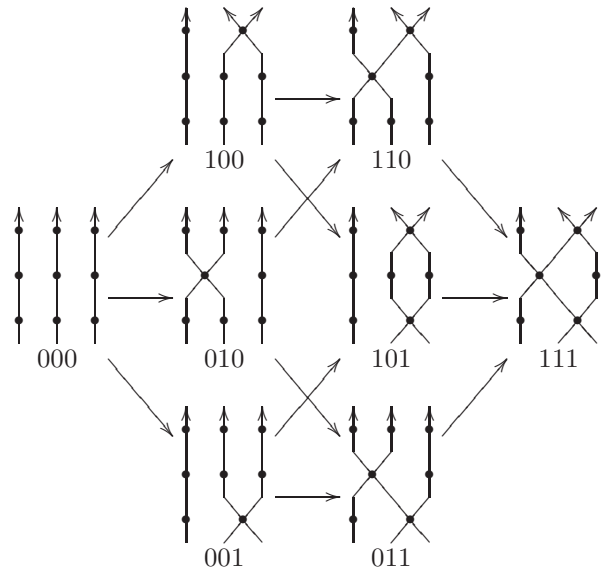}
\caption{Portion of the cube of resolutions for $D_2$ in layers $s_1$, $s_2$, and $s_3$.}
\label{reid3othercomplexesneg}
\end{center}
\end{figure}

Figures~\ref{reid3complexesneg} and~\ref{reid3othercomplexesneg} show the relevant portion of the cubes of resolutions associated to $D_1$ and $D_2$, respectively. Throughout this section, we will abbreviate indices to three places, relabeling the diagrams as necessary so that the Reidemeister III move occurs in layers $s_1$, $s_2$, and $s_3$, and using the index to indicate the states of the crossings in those layers only.


The goal is to prove that the chain complexes $C(D_1)$ and $C(D_2)$ are chain homotopy equivalent.  The strategy will be to prove that they are each chain homotopy equivalent to simplified complexes: $C(D_1)$ will be homotopic to

\begin{equation*}
\xymatrix@R-15pt{\, & \, & \mathcal{A}_{010}(D_1) \ar[r] \ar[ddr] & \mathcal{A}_{110}(D_1) \ar[dr] & \,\\
\overline{C}(D_1)= & \mathcal{A}_{000}(D_1) \ar[ur] \ar[dr] & \, & \, & \mathcal{B} \\
\, & \, & \mathcal{A}_{001}(D_1) \ar[uur] \ar[r] & \mathcal{A}_{011}(D_1) \ar[ur] & \, }
\end{equation*}

\noi and $C(D_2)$ will be homotopic to 

\begin{equation*}
\xymatrix@R-15pt{\, & \, & \mathcal{A}_{001}(D_2) \ar[r] \ar[ddr] & \mathcal{A}_{011}(D_2) \ar[dr] & \,\\
\overline{C}(D_2)= & \mathcal{A}_{000}(D_2) \ar[ur] \ar[dr] & \, & \, & \mathcal{B} \\
\, & \, & \mathcal{A}_{010}(D_2) \ar[uur] \ar[r] & \mathcal{A}_{110}(D_2) \ar[ur] & \, }
\end{equation*}

\noi and we will exhibit an isomorphism between the simplified complexes.

\begin{figure}[htbp]
\begin{center}
\input{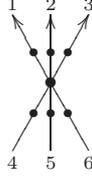}
\caption{The summand $\cB$ that appears in the simplifications of $C(D_1)$ and $C(D_2)$ is naturally associated to this diagram with a 6-valent vertex. See the beginning of Section~\ref{sec:reid3} and the last part in the proof of Lemma~\ref{reid3splitting111} for more details.}
\label{6valentvertex}
\end{center}
\end{figure}

The module $\mathcal{B}$ is a direct summand common to $\cA_{111}(D_1)$ and $\cA_{111}(D_2)$. It is naturally associated to the diagram $D_\ast$, in which the crossings in the Reidemeister III move and their nearby bivalent vertices are replaced by the 6-valent vertex shown in Figure~\ref{6valentvertex}. We define $\cB$ as a module over $\cR[x_0,\ldots,x_6,x_{13},\ldots,x_n]$ as follows. First, define the local relations associated to a 6-valent vertex using elementary symmetric polynomials:
\[\cL^{123}_{\text{sym}}=\left(\begin{array}{c}
t^3(x_1+x_2+x_3)-(x_4+x_5+x_6)\\
t^6\sigma_2(x_1,x_2,x_3)-\sigma_2(x_4,x_5,x_6)
\end{array}\right).\]
Next, define non-local relations for diagrams containing a 6-valent vertex using coherent regions, just as in Definition~\ref{dfnregions}, with the 6-valent vertex contributing a weight of 3 to the total weight of a coherent region. Let $\cN^{123}_{\text{sym}}$ be the ideal generated by non-local relations associated to coherent regions in $D_\ast$ that contain or have on their boundary the 6-valent vertex. Let $\cN^\prime$ be generated by non-local relations associated to all other coherent regions. Note that such regions and their associated relations are the same for $D_\ast$ as for $D_1$ and $D_2$ with any choice of resolution in layers $s_1$, $s_2$, and $s_3$. 

Define $\cB$ by
\[\cB=\frac{\cR[x_0,\ldots,x_6,x_{13},\ldots,x_n]}{\cL^\prime+\cL_{\text{sym}}^{123}+(t^9x_1x_2x_3-x_4x_5x_6)+\cN^\prime+\cN_{\text{sym}}^{123}},\]
where $\cL^\prime$ is generated by the local relations associated to layers $s_i$ for $i>3$. Note that $\cL^\prime$ is also the same for $D_\ast$ as for $D_1$ or $D_2$ with any choice of resolution in layers $s_1$, $s_2$, and $s_3$. 

The other modules in the simplified complexes $\overline{C}(D_i)$ correspond to resolutions of $D_1$ and $D_2$ that are identical after moving a layer of bivalent vertices. Specifically, the 010-resolution of $D_1$ matches the 001-resolution of $D_2$ and vice versa; and the 110-resolution of $D_1$ matches the 011-resolution of $D_2$ and vice versa. The 000-resolutions of $D_1$ and $D_2$ are identical.

The argument for simplifying $C(D_i)$ to $\overline{C}(D_i)$ proceeds in two parts. In Section~\ref{sec:reid3splittings}, we establish direct sum splittings of $\cA_{111}(D_i)$ and $\cA_{101}(D_i)$ for $i=1,2$, which categorify the MOY relations in the top and bottom lines of Figure~\ref{fig:braidmoy}, respectively. These splittings show that $C(D_1)$ and $C(D_2)$ both have the following form. 

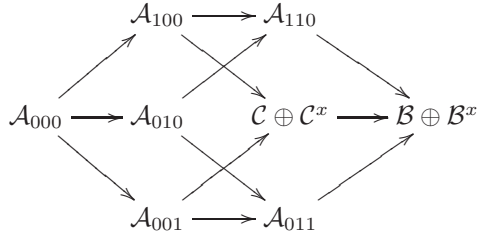
\begin{equation*}
\xymatrix@C-5pt{ \, & \mathcal{A}_{100} \ar[r] \ar[dr] & \mathcal{A}_{110} \ar[dr] & \,\\
\cA_{000} \ar[ur] \ar[r] \ar[dr] & \cA_{010} \ar[ur]_<<<<<<{\!\!\!} \ar[dr]^<<<<<<<{\!\!\!} & \mathcal{C}\oplus\mathcal{C}^x \ar[r] & \mathcal{B}\oplus\mathcal{B}^x \\
\, & \mathcal{A}_{001} \ar[ur] \ar[r] & \mathcal{A}_{011} \ar[ur] & \, }
\end{equation*}

After a series of changes of basis carried out in Section~\ref{sec:reid3complexmanipulation}, we identify contractible summands $\cA_{100}(D_i)\to\cC^x(D_i)$ and $\cC(D_i)\to\cB^x(D_i)$ in $C(D_i)$ for $i=1,2$. Removing these subcomplexes yields the simplified complexes $\overline{C}(D_i)$. Section~\ref{sec:reid3complexmanipulation} also exhibits an isomorphism between $\overline{C}(D_1)$ and $\overline{C}(D_2)$.

\subsection{Splittings for Reidemeister Move III}
\label{sec:reid3splittings}

This section establishes four direct sum splittings: two each for $D_1$ and $D_2$. In each case, the basic outline of the proof is similar for $D_1$ and $D_2$. We argue in full detail for the $D_1$ splittings, with periodic indications about how to modify these arguments for the $D_2$ cases. For arguments that apply equally well to $D_1$ and $D_2$, we typically omit $D_i$ from the notation (e.g.~using $\cA_{101}$ rather than $\cA_{101}(D_i)$).

\begin{lemma}
\label{reid3splitting101}
The algebras associated to the 101-resolutions of $D_1$ and $D_2$ split as direct sums of $\cR[x_0,\ldots,x_6,x_{13},\ldots,x_n]$-modules
 $\mathcal{A}_{101}(D_i)\isom\mathcal{C}(D_i)\oplus\mathcal{C}^x(D_i)$, where  $\mathcal{C}^x(D_i)\isom\mathcal{A}_{100}(D_i)$ and the edge map $\cA_{100}(D_i)\to\cA_{101}(D_i)$ is an isomorphism onto the $\cC^x(D_i)$ summand.
\end{lemma}

\begin{lemma}
\label{reid3splitting111}
The algebras associated to the 111-resolutions of $D_1$ and $D_2$ split as direct sums of $\cR[x_0,\ldots,x_6,x_{13},\ldots,x_n]$-modules
$\mathcal{A}_{111}(D_i)\isom\mathcal{B}\oplus\mathcal{B}^x(D_i)$, where $\cB$ is the module described above, $\mathcal{B}^x(D_i)\isom\mathcal{C}(D_i)$, and the edge map $\cA_{101}\to\cA_{111}$ restricted to the $\mathcal{C}(D_i)$ summand is an isomorphism onto $\cB^x(D_i)$.
\end{lemma}

\begin{proof}[Proof of Lemma~\ref{reid3splitting101}]
We know from Lemma~\ref{lemma:reid2} that $\mathcal{A}_{101}$ splits as a direct sum of modules isomorphic to $\mathcal{A}_{100}$ and $\mathcal{A}_{001}$.  However, it will be useful to establish a particular splitting so that we may see directly the isomorphisms $\mathcal{A}_{100}(D_i)\to \mathcal{C}^x(D_i)$ and (in the proof of Lemma~\ref{reid3splitting111}) $\mathcal{C}(D_i)\to\mathcal{B}^x(D_i)$.

The idea is similar to that of Lemma~\ref{lemma:reid2}.  We first manipulate the local relations from the vertices in $D_i$ to a convenient form, then obtain direct sum splittings by eliminating all quadratic and higher-order appearances of one variable, and keep track throughout of how these manipulations affect the non-local relations.

We begin with the presentation of $\cA_{101}$ as $$\cA_{101}\isom\frac{\cR[x_0,\ldots,x_n]}{\cL^{123}_{101}+\cL^\prime+\cN_{101}},$$ where $\cL_{101}^{123}$ is generated by local relations from layers $s_1$, $s_2$, and $s_3$, $\cL^\prime$ is generated by  local relations associated to other layers, and $\cN_{101}$ is generated by non-local relations.  Note that $\cL^\prime$ is generated by relations that do not use any of $x_7,\ldots,x_{12}$ and that it is identical for any resolution of layers $s_1$, $s_2$, and $s_3$ in either $D_1$ or $D_2$ and for $D_\ast$. It will not be affected by any of the calculations below. Thinking of non-local relations as coming from coherent regions, and applying Observation~\ref{regionobs}, notice that $x_7,\ldots,x_{12}$ need not ever appear in a generating set for $\mathcal{N}_{101}$. In $D_1$, any coherent region containing the elementary region to the right of $x_7$ and $x_9$ can be assumed to include the bigon bounded by edges $x_7$, $x_8$, $x_9$, and $x_{10}$. Similarly, in $D_2$, any coherent region containing the elementary region to the right of $x_8$ and $x_{10}$ may be assumed to include the bigon bounded by $x_7$, $x_8$, $x_9$, and $x_{10}$. The manipulations below will not affect such a generating set for $\cN_{101}$.

The following matrix is a generating set for $\cL^{123}_{101}(D_1)$, with $x_{11}$ and $x_{12}$ already eliminated using linear relations $tx_3-x_{11}$ and $tx_{12}-x_6$.  The analogous generating set for $\cL^{123}_{101}(D_2)$ can be produced by using $tx_1-x_{11}$ and $tx_{12}-x_4$ to eliminate $x_{11}$ and $x_{12}$. It is related to the matrix for $D_1$ below by exchanging $x_3$ with $x_1$ and $x_6$ with $x_4$ throughout.
$$
\left(\begin{array}{c}
t(x_1+x_2)-(x_7+x_8)\\
t^2x_1x_2-x_7x_8\\
t^3x_3-x_6\\
tx_7-x_9\\
t(x_9+x_{10})-(x_4+x_5)\\
(tx_9-x_5)(x_4-tx_9)\\
tx_8-x_{10}
\end{array}\right)$$
Use row IV to eliminate $x_7$, row VII to eliminate $x_8$, and row V to eliminate $x_{10}$, then rearrange. 
$$\xymatrix{
{\left(\begin{array}{c}
t(x_1+x_2)-t^{-2}(x_4+x_5)\\
t^2x_1x_2+t^{-2}x_9^2-t^{-3}x_9(x_4+x_5)\\
t^3x_3-x_6\\
(tx_9-x_5)(x_4-tx_9)\\
\end{array}\right)}
\ar[d]^{\text{I}+t^{-2}\text{III}\text{\,\,and\,\,}\text{II}+t^{-4}\text{IV}} 
\\
{\left(\begin{array}{c}
t(x_1+x_2+x_3)-t^{-2}(x_4+x_5+x_6)\\
t^2x_1x_2-t^{-4}x_4x_5\\
t^3x_3-x_6\\
(tx_9-x_5)(x_4-tx_9)\\
\end{array}\right)} 
}$$
Clear negative powers of $t$ from all rows and symmetrize the presentation as follows.
$$\xymatrix{
{\left(\begin{array}{c}
t^3(x_1+x_2+x_3)-(x_4+x_5+x_6)\\
t^6x_1x_2-x_4x_5\\
t^3x_3-x_6\\
(tx_9-x_5)(x_4-tx_9)\\
\end{array}\right)}
\ar[d]^{\text{II}+t^3(x_1+x_2)\text{III}+x_6\text{I}-x_6\text{III}}
\\
{\left(\begin{array}{c}
t^3(x_1+x_2+x_3)-(x_4+x_5+x_6)\\
t^6\sigma_2(x_1,x_2,x_3)-\sigma_2(x_4,x_5,x_6)\\
t^3x_3-x_6\\
(tx_9-x_5)(x_4-tx_9)\\
\end{array}\right)}
}$$
where $\sigma_2$ is the second elementary symmetric polynomial.

Apply the same sequence of operations to the presentation for $\cL_{101}^{123}(D_2)$, except replace $x_1$ with $x_3$ and $x_6$ with $x_4$ in the final row operation. The result is the following new presentation for $\cL_{101}^{123}(D_2)$.
$${\left(\begin{array}{c}
t^3(x_1+x_2+x_3)-(x_4+x_5+x_6)\\
t^6\sigma_2(x_1,x_2,x_3)-\sigma_2(x_4,x_5,x_6)\\
t^3x_1-x_4\\
(tx_9-x_5)(x_6-tx_9)\\
\end{array}\right)}
$$

Notice that the first two rows in the new presentations of $\cL_{101}^{123}(D_1)$ and $\cL_{101}^{123}(D_2)$ above are familiar as the generators of $\cL^{123}_{\text{sym}}$ in the definition of $\cB$. Notice also that $\cL^{123}_{\textrm{sym}}$ is generated by relations that do not use any of $x_7,\ldots,x_{12}$. Define $$\cT=\frac{\mathcal{R}[x_0,\ldots,x_6,x_{13},\ldots,x_n]}{\cL^\prime+\cL^{123}_{\textrm{sym}}}.$$ This definition works equally well for $D_1$, $D_2$, and $D_\ast$.

Let $q_1=(tx_9-x_5)(x_4-tx_9)$ and $q_2=(tx_9-x_5)(x_6-tx_9)$. These are the only relations in the presentations of $\cL_{101}^{123}(D_i)$ above that use any of $x_7,\ldots,x_{12}$. Let $r_1=t^3x_3-x_6$ and $r_2=t^3x_1-x_4$. So far, we have established that $$\cA_{101}(D_i)\isom\frac{\cT[x_9]}{(q_i)+ (r_i)+\cN_{101}(D_i)}$$ and that $x_9$ appears only in $q_i$.  Since $q_i$ is quadratic in $x_9$, we could use it to replace any appearance of $x_9^k$ for $k\geq 2$ in a presentation of $\cA_{101}(D_i)$ with some polynomial that was linear in $x_9$. However, we have already eliminated all appearances of $x_9$ from the rest of the presentation.  Therefore, we may instead forget the relation $q_i$, and split $\mathcal{A}_{101}(D_i)$ into a summand generated by $1$ and a summand generated by a polynomial that is linear in $x_9$.
\begin{align*}
\mathcal{A}_{101}(D_1)&\isom\frac{\cT(1)}{(r_1)+\cN_{101}(D_1)}\bigoplus\frac{\cT(tx_9-x_4)}{(r_1)+\cN_{101}(D_1)}\\
\mathcal{A}_{101}(D_2)&\isom\frac{\cT(1)}{(r_2)+\cN_{101}(D_2)}\bigoplus\frac{\cT(tx_9-x_6)}{(r_2)+\cN_{101}(D_2)}.
\end{align*}
With the first summand in each case as $\mathcal{C}(D_i)$ and the second as $\mathcal{C}^x(D_i)$, this is the splitting asserted in the statement of the lemma.

We now check that $\mathcal{A}_{100}\isom\mathcal{C}^x$ via the edge map $\cA_{100}\to\cA_{101}$, which is multiplication by $tx_9-x_4$ for $D_1$ and by $tx_9-x_6$ for $D_2$. The edge map definitely takes the generator ($1$) of $\cA_{100}$ to the generator of the $\cC^x$ summand of $\cA_{101}$. To check that this is an isomorphism of $\cT$-modules (hence of $\cR[x_0,\ldots,x_6,x_{13},\ldots,x_n]$-modules), we simplify the presentation of $\cA_{100}$ and match it with that of $\cC^x$. 

Begin with the presentation of $\cA_{100}$ as 
\[\cA_{100}\isom\frac{\cR[x_0,\ldots,x_n]}{\cL_{100}^{123}+\cL^\prime+\cN_{100}},\] with notation for the ideals analogous to that used in the presentation of $\cA_{101}$. Note that $x_7,\ldots,x_{12}$ do not appear in the standard generating set for $\cL^\prime$, since it concerns only layers $s_i$ for $i>3$. These variables also need not appear in a minimal generating set for $\mathcal{N}_{100}$. If a subset had one of these as an outgoing or incoming edge, we could use the relations associated to bivalent vertices to eliminate it. 

Turning to $\cL^{123}_{100}(D_1)$, the local relations from layers $s_1$, $s_2$, and $s_3$, eliminate $x_{11}$ and $x_{12}$ immediately using the linear relations on the rightmost strand, then remove $x_7,\ldots,x_{10}$ as follows.
$$\xymatrix{
{\left(\begin{array}{c}
tx_1+tx_2-x_7-x_8\\
(tx_2-x_7)(x_8-tx_2)\\
t^3x_3-x_6\\
tx_7-x_9\\
tx_8-x_{10}\\
tx_9-x_5\\
tx_{10}-x_4
\end{array}\right)}
\ar[d]^{\text{I}+t^{-2}\text{III}+t^{-1}\text{IV}+t^{-1}\text{V}+t^{-2}\text{VI}+t^{-2}\text{VII}}
\\
{\left(\begin{array}{c}
t(x_1+x_2+x_3)-t^{-2}(x_4+x_5+x_6)\\
(tx_2-x_7)(x_8-tx_2)\\
t^3x_3-x_6\\
tx_7-x_9\\
tx_8-x_{10}\\
tx_9-x_5\\
tx_{10}-x_4
\end{array}\right)}
}$$
\noindent Simplify by multiplying the first row by $t^2$, using row $\text{IV}$ to eliminate $x_7$, row $\text{V}$ to eliminate $x_8$, row $\text{VI}$ to eliminate $x_9$, and row $\text{VII}$ to eliminate $x_{10}$. Then multiply the second row by $t^4$.
$$\xymatrix{
{\left(\begin{array}{c}
t^3(x_1+x_2+x_3)-(x_4+x_5+x_6)\\
(t^3x_2-x_5)(x_4-t^3x_2)\\
t^3x_3-x_6
\end{array}\right)}
\ar[d]^{\text{II}+t^3(x_2+x_3)\text{I}+(x_4+x_5-t^3(x_2+x_3))\text{III}}
\\
{\left(\begin{array}{c}
t^3(x_1+x_2+x_3)-(x_4+x_5+x_6)\\
t^6\sigma_2(x_1,x_2,x_3)-\sigma_2(x_4,x_5,x_6)\\
t^3x_3-x_6
\end{array}\right)}
}$$

A similar computation in the $D_2$ case produces the following presentation for $\cL_{100}^{123}(D_2)$
\[{\left(\begin{array}{c}
t^3(x_1+x_2+x_3)-(x_4+x_5+x_6)\\
t^6\sigma_2(x_1,x_2,x_3)-\sigma_2(x_4,x_5,x_6)\\
t^3x_1-x_4
\end{array}\right)}
\]

In each case, the top two rows generate $\cL^{123}_{\textrm{sym}}$ and the bottom row is $r_i$ as previously defined. Therefore, we have established that $$\cA_{100}(D_i)\isom\frac{\cT}{(r_i)+\cN_{100}(D_i)}.$$  It remains to check that $\mathcal{N}_{100}=\mathcal{N}_{101}$.  Figure~\ref{100v101cycles} shows how the cycles that pass through the 101-resolution of $D_1$ pair up with the cycles that pass through the 100-resolution of $D_1$ to give equivalent non-local relations. There is an analogous way of pairing cycles in the 101- and 100-resolutions of $D_2$. Any cycle that does not pass through this region certainly has the same associated non-local relation in $\mathcal{N}_{101}$ and $\mathcal{N}_{100}$.  We have identified identical generating sets for $\mathcal{N}_{101}$ and $\mathcal{N}_{100}$. Therefore, $\cA_{100}$ and $\cC^x$ have identical presentations as $\cT$-modules. As previously noticed, the edge map from $\cA_{100}$ to $\cA_{101}$ sends the generator $(1)$ of $\cA_{100}$ to the generator $(tx_9-x_4)$ or $(tx_9-x_6)$ of $\cC^x(D_i)$, as appropriate. We conclude that this edge map is an isomorphism onto $\cC^x(D_i)$.
\begin{figure}[htbp]
\begin{center}
\input{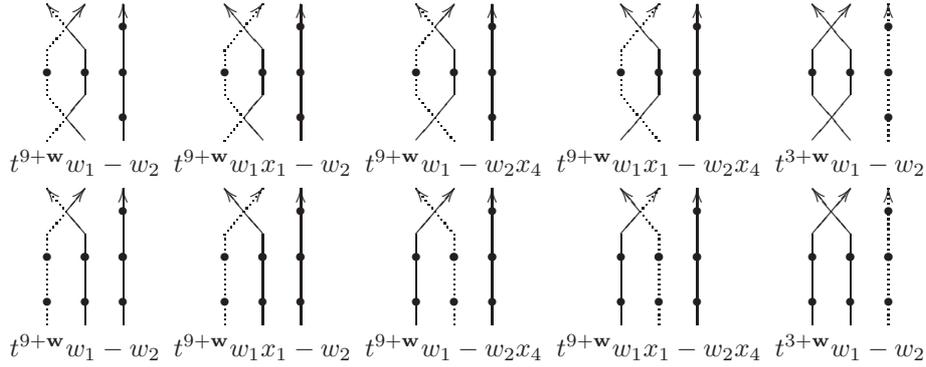}
\caption{Pairing of cycles (dotted lines) that pass through the 101-resolution (top row) and the 100-resolution (bottom row) of $D_1$. Assume that the braid axis is to the right of each picture. We use $\weight $ to denote the weight, $w_1$ the product of outgoing edges, and $w_2$ the product of incoming edges for the portion of the cycle away from the Reidemeister III move. An analogous matching can be made for cycles passing through the 101- and 100-resolutions of $D_2$.}
\label{100v101cycles}
\end{center}
\end{figure}
\end{proof}

The proof of Lemma~\ref{reid3splitting111} is similar, except that more work is required to keep track of the non-local relations.  As before, we use local relations associated to layers $s_1$, $s_2$, and $s_3$ in $\mathcal{A}_{111}$ to eliminate several edge variables, then use a quadratic relation to split $\cA_{111}$ as a direct sum, and finally check that one of the direct summands is isomorphic to $\mathcal{C}$ via the appropriate edge map while the other is isomorphic to $\cB$. As before, we omit $D_i$ from the notation when an argument applies to both $D_1$ and $D_2$.

\begin{proof}[Proof of Lemma~\ref{reid3splitting111}]
Let $\mathcal{L}^{123}_{111}$ denote the ideal generated by local relations associated to layers $s_1$, $s_2$, and $s_3$, while $\mathcal{L}^\prime$ denotes the ideal generated by local relations associated to all other layers as before. Let $\cN_{111}$ denote the ideal generated by non-local relations in the 111-resolution. So we begin with $$\cA_{111}\isom\frac{\cR[x_0,\ldots,x_n]}{\mathcal{L}^{123}_{111}+\mathcal{L}^\prime+\cN_{111}}.$$

The general strategy will be to eliminate $x_7$, $x_8$, $x_{10}$, $x_{11}$, and $x_{12}$ from the presentation of $\cA_{111}$ and limit use of $x_9$ as much as possible. We will then rewrite $\cA_{111}$ in the form $\cT[x_9]/\mathcal{I}$ for an appropriate ideal $\mathcal{I}$, where $\cT$ is the $\cR[x_0,\ldots,x_6,x_{13},\ldots,x_n]$-algebra defined in the proof of Lemma~\ref{reid3splitting101}. Finally, we will use the quadratic relation associated to layer $s_3$, denoted $q_i$ as in the proof of Lemma~\ref{reid3splitting101}, to split $\cA_{111}$ into direct summands generated by $1$ and a linear polynomial in $x_9$.

Notice first that no part of this strategy will affect the ideal $\mathcal{L}^\prime$. Edges $x_7,\ldots,x_{12}$ connect layer $s_1$ to layer $s_2$ or layer $s_2$ to layer $s_3$, so they do not appear in local relations associated to any other layers.

\subsection*{Analysis of Local Relations}
\noi For the presentation of $\mathcal{L}^{123}_{111}$, first use relations associated to bivalent vertices to replace $x_{11}$ and $x_{12}$. Then the following matrix is a generating set for $\mathcal{L}^{123}_{111}(D_1)$. The analogous generating set for $\cL^{123}_{111}(D_2)$ is related to the matrix below by exchanging $x_1$ with $x_3$ and $x_4$ with $x_6$.
$$\xymatrix@R-15pt{
{\left(\begin{array}{c}
t(x_1+x_2)-(x_7+x_8)\\
(tx_2-x_7)(x_8-tx_2)\\
t(tx_3+x_7)-(t^{-1}x_6+x_9)\\
(t^3x_3-x_6)(t^{-1}x_9-tx_3)\\
t(x_9+x_{10})-(x_4+x_5)\\
(tx_9-x_5)(tx_9-x_4)\\
tx_8-x_{10}
\end{array}\right)}
\ar[d]^{\text{I}+t^{-1}\text{III}+t^{-2}\text{V}+t^{-1}\text{VII}}\\
{\left(\begin{array}{c}
t(x_1+x_2+x_3)-t^{-2}(x_4+x_5+x_6)\\
(tx_2-x_7)(x_8-tx_2)\\
t(tx_3+x_7)-(t^{-1}x_6+x_9)\\
(t^3x_3-x_6)(t^{-1}x_9-tx_3)\\
t(x_9+x_{10})-(x_4+x_5)\\
(tx_9-x_5)(tx_9-x_4)\\
tx_8-x_{10}
\end{array}\right)}
}$$

\noi Multiply row I by $t^2$, use row III to eliminate $x_7$, use row V to eliminate $x_{10}$, and use row VII to eliminate $x_8$, then multiply row II by $t^4$.
$$
\left(\begin{array}{c}
t^3(x_1+x_2+x_3)-(x_4+x_5+x_6)\\
(t^3(x_2+x_3)-x_6-tx_9)(x_4+x_5-tx_9-t^3x_2)\\
(t^3x_3-x_6)(t^{-1}x_9-tx_3)\\
(tx_9-x_5)(tx_9-x_4)\\
\end{array}\right)$$
Use row IV to replace $t^2x_9^2$ in row II, then add $t^2\text{III}$ and $t^3(x_2+x_3)\text{I}$ to row II, and then multiply row III by $t$ to obtain
$$
\left(\begin{array}{c}
t^3(x_1+x_2+x_3)-(x_4+x_5+x_6)\\
t^6\sigma_2(x_1,x_2,x_3)-\sigma_2(x_4,x_5,x_6)\\
(t^3x_3-x_6)(x_9-t^2x_3)\\
(tx_9-x_5)(tx_9-x_4)\\
\end{array}\right)
$$
where $\sigma_2$ is the second elementary symmetric polynomial. 

We may apply an analogous sequence of transformations to the presentation of $\cL_{111}^{123}(D_2)$, at each step maintaining the property that the presentations of $\cL_{111}^{123}(D_i)$ are related by exchanging $x_1$ with $x_3$ and $x_4$ with $x_6$. This analogous sequence of transformations is in fact identical until the last step, in which we should add $t^2\text{III}$ and $t^2(x_1+x_2)\text{I}$ to row $\text{II}$. To be explicit, we obtain the following presentation for $\cL_{111}^{123}(D_2)$.
\[\left(\begin{array}{c}
t^3(x_1+x_2+x_3)-(x_4+x_5+x_6)\\
t^6\sigma_2(x_1,x_2,x_3)-\sigma_2(x_4,x_5,x_6)\\
(t^3x_1-x_4)(x_9-t^2x_1)\\
(tx_9-x_5)(tx_9-x_6)\\
\end{array}\right)
\]

In both presentations, the first two relations generate the familiar ideal $\mathcal{L}^{123}_{\textrm{sym}}$. The last relation in both presentations is familiar as $q_i$ from the proof of Lemma~\ref{reid3splitting101}. Let $p_1=(t^3x_3-x_6)(x_9-t^2x_3)$ and $p_2=(t^3x_1-x_4)(x_9-t^2x_1)$. Then we have expressed $\mathcal{L}^{123}_{111}$ as the sum of an ideal whose generating set does not involve $x_9$, the quadratic relation $q_i$ that will be used to split $\cA_{111}$ as a direct sum, and the relation $p_i$, which we will have to follow up carefully. We may retain the definition of $\cT$ from the proof of Lemma~\ref{reid3splitting101}, and write 
$$\cA_{111}(D_i)\isom\frac{\cT[x_7,\ldots,x_{12}]}{(p_i)+(q_i)+\cN_{111}(D_i)}.$$

\subsection*{Analysis of Non-local Relations: The $D_1$ Case}
\noi We turn next to an examination of the ideal $\cN_{111}$ of non-local relations, which will require separate arguments for $D_1$ and $D_2$. Label the elementary regions in the vicinity of the Reidemeister III move as in Figure~\ref{fig:111elemregions}. As usual, assume that the braid axis is to the right of each diagram. 

\begin{figure}[tbp]
\begin{center}
\input{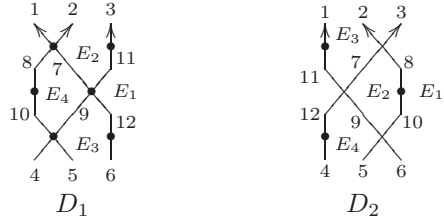}
\caption{Elementary regions in the vicinity of the Reidemeister III move in the 111-resolution of $D_1$ (left) and $D_2$ (right).}
\label{fig:111elemregions}
\end{center}
\end{figure}

Using Definition~\ref{dfnregions} for the generators of $\cN_{111}(D_1)$, we can split $\cN_{111}(D_1)$ into a sum of five ideals based on types of coherent regions.  Let $\cN^\prime$ be the ideal generated by the relations from coherent regions that do not use any of $E_1$, $E_2$, $E_3$, or $E_4$. Notice that this ideal is identical for $D_1$ and $D_2$, even though $E_i$ actually refers to different elementary regions in each diagram. Moreover, none of the relations associated to such coherent regions use edge variables $x_7,\ldots,x_{12}$, so they will carry through all of our calculations unchanged. 

Let $\cE^{1234}(D_1)$ be generated by relations from coherent regions that use all of $E_1$, $E_2$, $E_3$, and $E_4$. These relations may involve $x_1$ and $x_4$, but not any of $x_7,\ldots,x_{12}$, so they carry through our calculations unchanged. The ideal $\cE^{1234}(D_1)$ also accounts for relations associated to coherent regions that contain $E_1$, $E_2,$ and $E_3$. Adding $E_4$ to such a region would add only the bivalent vertex between edges 8 and 10, which is exactly the situation described in Observation~\ref{regionobs}. Therefore, we need not consider coherent regions that contain $E_1$, $E_2,$ and $E_3$ without $E_4$ in a minimal generating set for $\cN_{111}(D_1)$.

Let $\cE^{12}(D_1)$ (respectively $\cE^{13}(D_1)$) be generated by non-local relations from coherent regions that use $E_1$ and $E_2$, but not $E_3$ or $E_4$ (respectively $E_1$ and $E_3$ but not $E_2$ or $E_4$). Some of the edge variables $x_7,\ldots,x_{12}$ do appear in the relations associated to such regions, but can be easily eliminated using the quadratic relation from layers $s_1$ or $s_3$ as appropriate.  Figure~\ref{111cycles1213} shows the necessary calculations in each case.

\begin{figure}[tbp]
\begin{center}
\input{111cycles1213}
\caption{Relations that generate $\cE^{12}(D_1)$ and $\cE^{13}(D_1)$, along with modifications to avoid the use of $x_7,\ldots,x_{12}$.  Assume that the braid axis is to the right of each picture. Brackets around $x_i$ denote an edge variable that may or may not occur in a relation depending on whether the coherent region under consideration contains an elementary region immediately above or below $E_4$. Dotted lines show the boundary of the relevant coherent region when such adjacent elementary regions are not included. In each diagram, $\weight $, $w_{\text{out}}$, and $w_{\text{in}}$ come from the portions of the region not shown in these local pictures.}
\label{111cycles1213}
\end{center}
\end{figure}

Finally, let $\cE^1(D_1)$ be generated by relations from coherent regions that use $E_1$ but none of $E_2$, $E_3$, or $E_4$, as shown in Figure~\ref{111cycles1}. These relations have the form $t^{4+\weight }w_{\text{out}}x_7-w_{\text{in}}x_9$, where $\weight $, $w_{\text{out}}$, and $w_{\text{in}}$ come from the portions of the coherent region not shown in Figure~\ref{111cycles1}. We will not be able to simultaneously eliminate $x_7,\ldots,x_{12}$ from such a relation, but we can eliminate all except $x_9$ using the linear relations from the crossing in layer $s_2$ and linear relations associated to bivalent vertices. In fact, we can write any generator of $\cE^1(D_1)$ in the form $t^{2+\weight }w_{\text{out}}(x_6-t^3x_3)+x_9(t^{3+\weight }w_{\text{out}}-w_{\text{in}})$, where $w_{\text{out}}$ and $w_{\text{in}}$ are words in $x_0,\ldots,x_6,x_{13},\ldots,x_n$ only.  

\begin{figure}[tbp]
\begin{center}
\input{111cycles1}
\caption{Removing $x_7$ from relations that generate $\cE^1(D_1)$.}
\label{111cycles1}
\end{center}
\end{figure}

We have exhausted the possible combinations of elementary regions $E_1,\ldots,E_4$ that can appear in a coherent region, so we may now express $\cN_{111}(D_1)$ as $$\cN_{111}(D_1)=\cN^\prime+\cE^{1234}(D_1)+\cE^{12}(D_1)+\cE^{13}(D_1)+\cE^1(D_1).$$ Moreover, we have eliminated all appearances of $x_7,\ldots,x_{12}$ from the generating sets of $\cN^\prime$, $\cE^{1234}(D_1)$, $\cE^{12}(D_1),$ and $\cE^{13}(D_1).$ Defining $\cT^\prime_1$ by $$\cT^\prime_1=\frac{\cT}{\cN^\prime+\cE^{12}(D_1)+\cE^{13}(D_1)+\cE^{1234}(D_1)}$$ we then have a presentation of $\cA_{111}(D_1)$ as $$\cA_{111}(D_1)\isom\frac{\cT^\prime_1[x_9]}{(p_1)+(q_1)+\cE^1(D_1)}.$$  

The next step will be to use $q_1$ to split $\cA_{111}(D_1)$ as a direct sum of modules over $\cR[x_0,\ldots,x_6,x_{13},\ldots,x_n]$, one of which is generated by (1) and one of which is generated by $(t^2x_3-x_9)$. In other words, we would like to find ideals $\cP(D_1),\cP^x(D_1),\cE(D_1),$ and $\cE^x(D_1)$ in $\cT^\prime_1$ such that 
 $$\frac{\cT^\prime_1[x_9]}{(p_1)+(q_1)+\cE^1(D_1)}\isom\frac{\cT^\prime_1(1)}{\cP(D_1)+\cE(D_1)}\bigoplus\frac{\cT^\prime_1(t^2x_3-x_9)}{\cP^x(D_1)+\cE^x(D_1)}$$ as $\cR[x_0,\ldots,x_6,x_{13},\ldots,x_n]$-modules.

As in the proof of Lemma~\ref{reid3splitting101}, we may use $q_i$ to replace any appearance of $x_9^k$ for $k\geq 2$ with a polynomial that is linear in $x_9$. This procedure has no effect on the ideals from which $x_9$ has been eliminated, but it does affect $(p_1)$ and $\cE^1(D_1)$.  To analyze how, think of the ideal that $p_1$ generates in $\cT^\prime_1[x_9]/(q_1)$ as the sum of the ideals generated by $p_1$ and $x_9p_1$. If we use $q_1$ to eliminate any appearances of $x_9^2$ in these generating sets, then we can find appropriate generators for $\cP(D_1)$ and $\cP^x(D_1)$ by writing $p_1$ and $x_9p_1$ in terms of 1 and $t^2x_3-x_9$. Actually, $p_1=(t^3x_3-x_6)(x_9-t^2x_3)$ is already in the correct format, so let $t^3x_3-x_6$ be one of the generators of $\cP^x(D_1)$. For $x_9p_1$, we calculate as follows, replacing $x_9^2$ using $q_1$, then eliminating a term using $p_1$.
\begin{align}
&\,\,x_9(x_9-t^2x_3)(t^3x_3-x_6)\nonumber\\
=&\,\,(t^{-1}x_9x_4+t^{-1}x_9x_5-t^{-2}x_4x_5-t^2x_3x_9)(t^3x_3-x_6)&\nonumber\\
=&\,\,(x_9-t^2x_3)(t^{-1}x_4+t^{-1}x_5-t^2x_3)(t^3x_3-x_6)\nonumber\\
&\,-\,(t^{-2}x_4x_5+t^4x_3^2-tx_3x_4-tx_3x_5)(t^3x_3-x_6)\nonumber\\
\equiv&\,\,(t^6x_3^2+x_4x_5-t^3x_3x_4-t^3x_3x_5)(t^3x_3-x_6)\nonumber\\
=&\,\,(t^3x_3-x_4)(t^3x_3-x_5)(t^3x_3-x_6)\label{square}
\end{align}
Therefore, the ideal generated by $p_1$ and $x_9p_1$ in $\cT^{\prime}_1[x_9]/(q_1)$ is equal to the ideal generated by $p_1$ and $(t^3x_3-x_4)(t^3x_3-x_5)(t^3x_3-x_6)$, which no longer uses $x_9$. It will be convenient to express this relation more symmetrically by modifying it using generators of $\cL^{123}_{\textrm{sym}}$. (Recall that the original $\cT$, of which $\cT^\prime_1$ is a quotient, was a quotient by $\cL^{123}_{\textrm{sym}}$, among other ideals.)
\begin{align*}
&(t^3x_3-x_4)(t^3x_3-x_5)(t^3x_3-x_6)\\ 
\equiv\,&(t^3x_3-x_4)(t^3x_3-x_5)(t^3x_3-x_6)\\
&+\, t^3x_3\left(t^6\sigma_2(x_1,x_2,x_3)-\sigma_2(x_4,x_5,x_6)\right)\\
&-\,  t^6x_3^2\left(t^3(x_1+x_2+x_3)-(x_4+x_5+x_6)\right)\\
=\,&t^9x_1x_2x_3-x_4x_5x_6
\end{align*}
 Let $\cP(D_1)$ be the ideal generated by $t^9x_1x_2x_3-x_4x_5x_6$ in $\cT^\prime_1$.  Adding the generator of $\cP(D_1)$ to $\cP^x(D_1)$ would not change the ideal, since $\cP^x(D_1)$ already has $t^3x_3-x_6$ as a generator.

We use the same strategy to find appropriate generators for $\cE(D_1)$ and $\cE^x(D_1)$. Generators of $\cE^1(D_1)$ have the form $f_1=t^{4+\weight}w_{\text{out}}x_7-w_{\text{in}}x_9$. We would like to write $f_1$ and $x_9f_1$ in terms of 1 and $t^2x_3-x_9$. We have already seen that $$f_1\equiv t^{2+\weight }w_{\text{out}}\left(x_6-t^3x_3\right)+x_9\left(t^{3+\weight }w_{\text{out}}-w_{\text{in}}\right),$$ where $w_{\text{out}}$ and $w_{\text{in}}$ are words in $x_0,\ldots,x_6,x_{13},\ldots,x_n$ only. Factoring out $x_9-t^2x_3$ yields $$f_1\equiv  (x_9-t^2x_3)(t^{3+\weight }w_{\text{out}}-w_{\text{in}})+t^2(t^\weight w_{\text{out}}x_6-w_{\text{in}}x_3).$$ 

Conveniently, the second term is a multiple of a generator of $\cN^\prime$ obtained as follows. Suppose $f_1$ came from a coherent region $R$. Let $V_R$ be the set of vertices contained in the closure of $R$, so that $f_1$ is the relation associated to $V_R$ under the subset interpretation of the non-local relations. Delete from $V_R$ the 4-valent vertex in layer $s_2$, the bivalent vertex between edges 3 and 11, and the bivalent vertex between edges 12 and 6. These deletions drop the weight of $V_R$ by 4. The resulting set of vertices has the same incoming and outgoing edges as $V_R$ except that $x_7$ has been replaced by $x_6$ and $x_9$ has been replaced by $x_3$. Therefore, the relation associated to this subset, which must appear in $\cN^\prime$, is exactly $t^\weight w_{\text{out}}x_6-w_{\text{in}}x_3$. As an element of $\cT^\prime_1$, the above expression for $f_1$ then simplifies to 
\begin{equation}
\label{finalnonlocal}
f_1\equiv\,(x_9-t^2x_3)(t^{3+\weight }w_{\text{out}}-w_{\text{in}}).
\end{equation} We conclude that a generating set for $\cE^x(D_1)$ should include $t^{3+\weight }w_{\text{out}}-w_{\text{in}}$. 

Next consider $x_9f_1$, using the expression for $f_1$ obtained in Equation~\ref{finalnonlocal} and the expression for $x_9(x_9-t^2x_3)$ obtained in the third line of the calculation preceding Equation~\ref{square}.
\begin{align}
\label{nonlocalsquare}
x_9f_1 =& \,x_9(x_9-t^2x_3)(t^{3+\weight }w_{\text{out}}-w_{\text{in}})\nonumber\\
=&\, (x_9-t^2x_3)(t^{-1}x_4+t^{-1}x_5-t^2x_3)(t^{3+\weight}w_{\text{out}}-w_{\text{in}})\nonumber\\
&-\,t^{-2}(t^3x_3-x_4)(t^3x_3-x_5)(t^{3+\weight}w_{\text{out}}-w_{\text{in}})\nonumber\\
\equiv& \,(t^3x_3-x_4)(t^3x_3-x_5)(t^{3+\weight }w_{\text{out}}-w_{\text{in}})
\end{align}
The last equivalence follows because the term we have eliminated is a multiple of $f_1$ as expressed in Equation~\ref{finalnonlocal}. These calculations eliminate all appearances of $x_7,\ldots,x_{12}$ from $x_9f_1$. We may again use relations in $\cL^{123}_{\text{sym}}$ to rewrite this expression in a more convenient form.
\begin{align*}
&\,(t^3x_3-x_4)(t^3x_3-x_5)(t^{3+\weight }w_{\text{out}}-w_{\text{in}})\\
\equiv& \,(t^3x_3-x_4)(t^3x_3-x_5)(t^{3+\weight }w_{\text{out}}-w_{\text{in}})\\
&-\,t^{6+\weight}w_{\text{out}}x_3\left(t^3x_1+t^3x_2+t^3x_3-x_4-x_5-x_6\right)\\
&+\,t^{3+\weight}w_{\text{out}}\left(t^6x_1x_2+t^6x_1x_3+t^6x_2x_3-x_4x_5-x_4x_6-x_5x_6\right)\\
=&\,t^3(t^3x_3-x_4-x_5)(t^{\weight}w_{\text{out}}x_6-w_{\text{in}}x_3)+(t^{9+\weight}w_{\text{out}}x_1x_2-w_{\text{in}}x_4x_5)
\end{align*}
The final expression above has $x_9f_1$ as a linear combination of two relations that are actually already accounted for in the splitting. The first relation, $t^{\weight}w_{\text{out}}x_6-w_{\text{in}}x_3$, arose during our analysis of $f_1$ above and was seen to be contained in $\cN^\prime$. The second, $t^{9+\weight}w_{\text{out}}x_1x_2-w_{\text{in}}x_4x_5$, is associated to a disconnected subset of vertices. If $f$ was associated to a region $R$ with corresponding subset of vertices $V_R$, then the second relation above is associated to the union of $V_R$ with the bivalent vertex between $x_8$ and $x_{10}$. Since that union is a disconnected subset, its associated relation is not needed in a minimal presentation of $\cN_{111}(D_1)$ (as noted just after Definition~\ref{dfnsubsets} in Section~\ref{sec:nonlocalrelations}). Therefore, it turns out that $x_9f_1$ does not add any new generators to $\cE(D_1)$ or $\cE^x(D_1)$. In fact, we did not need to put any generators in $\cE(D_1)$ at all.

We have now split $\cA_{111}(D_1)$ as a direct sum of $\cR[x_0,\ldots,x_6,x_{13},\ldots,x_n]$-modules:
$$\cA_{111}(D_1)\isom\frac{\cT^\prime_1(1)}{\mathcal{P}}\bigoplus\frac{\cT^\prime_1(t^2x_3-x_9)}{\mathcal{P}^x(D_1)+\cE^x(D_1)}.$$ Define $\cB(D_1)$ to be the first summand and $\cB^x(D_1)$ to be the second.

\subsection*{Analysis of Non-local Relations: The $D_2$ Case}
\noi Next, we analyze the non-local relations in the 111-resolution of $D_2$. Referring back to the labels in Figure~\ref{fig:111elemregions}, let $\cE^1(D_2)$, $\cE^{12}(D_2)$, $\cE^{123}(D_2)$, $\cE^{124}(D_2)$, and $\cE^{1234}(D_2)$ denote the ideals generated by non-local relations associated to coherent regions using the super-scripted elementary regions. Coherent regions cannot contain other combinations of elementary regions $E_1,\ldots,E_4$, so \[\cN_{111}(D_2)=\cN^\prime+\cE^1(D_2)+\cE^{12}(D_2)+\cE^{123}(D_1)+\cE^{124}(D_2)+\cE^{1234}(D_2).\] For all of these ideals except $\cE^1(D_2)$, there is a generating set that does not use any of $x_7,\ldots,x_{12}$. Simply use the relations $tx_1-x_{11}$ and $tx_{12}-x_4$ associated to bivalent vertices to eliminate $x_{11}$ and $x_{12}$ whenever they appear in one of the standard generators. 

As was the case for $D_1$, we cannot simultaneously eliminate all of $x_7,\ldots,x_{12}$ from a generating set of $\cE^1(D_2)$, but we can eliminate all except $x_9$. The standard generator in $\cE^1(D_2)$ has the form $f_2=t^{5+\weight}w_{\text{out}}x_9[x_2]-w_{\text{in}}x_7[x_5]$, where $[x_i]$ indicates edge variables that may or may not appear depending on whether the coherent region under consideration contains the elementary regions immediately above and below $E_2$. Replace $x_7$ using the linear relation $t(tx_1+x_7)-(t^{-1}x_4+x_9)$, which comes from the vertex in layer $s_2$ together with the bivalent vertices on the leftmost strand. Then regroup and clear negative powers of $t$ to obtain 
\begin{align}
\label{nonlocal2}
f_2&=t^{5+\weight}w_{\text{out}}x_9[x_2]-w_{\text{in}}x_7[x_5]\nonumber\\
&\equiv\,(t^3x_1-x_4)w_{\text{in}}[x_5]+tx_9(t^{6+\weight}w_{\text{out}}[x_2]-w_{\text{in}}[x_5]).\nonumber
\end{align}

Define $\cT^\prime_2$ to be the $\cT$-module
\[\cT^\prime_2=\frac{\cT}{\cN^\prime+\cE^{12}(D_2)+\cE^{123}(D_2)+\cE^{124}(D_2)+\cE^{1234}(D_2)}.\]
Then we have a presentation of $\cA_{111}(D_2)$ as 
\[\cA_{111}(D_2)\isom\frac{\cT^\prime_2[x_9]}{(p_2)+(q_2)+\cE^1(D_2)}.\]
As in the $D_1$ case, we would now like to find ideals $\cP(D_2)$, $\cP^x(D_2)$, $\cE(D_2)$, and $\cE^x(D_2)$ in $\cT^\prime_2$ such that 
\[\frac{\cT^\prime_2[x_9]}{(p_2)+(q_2)+\cE^1(D_2)}\isom\frac{\cT^\prime_2(1)}{\cP(D_2)+\cE(D_2)}\bigoplus\frac{\cT^\prime_2(t^2x_1-x_9)}{\cP^x(D_2)+\cE^x(D_2)}\] as $\cR[x_0,\ldots,x_6,x_{13},\ldots,x_n]$-modules.

The analysis is analogous to that of the $D_1$ case. We examine the polynomials $p_2$, $x_9p_2$, $f_2$, and $x_9f_2$, using $q_2$ to replace any appearance of $x_9^2$ and writing each polynomial in terms of $1$ and $t^2x_1-x_9$. We already have $p_2=(t^3x_1-x_4)(x_9-t^2x_1)$ in the correct format, so we add $t^3x_1-x_4$ as a generator of $\cP^x(D_2)$. It turns out that \[x_9p_2\equiv (t^3x_1-x_4)(t^3x_1-x_5)(t^3x_1-x_6)\] in $\cT^\prime_2[x_9]/(q_2)$. A quick modification by relations in $\cL^{123}_{\textrm{sym}}$ recovers the same nicely symmetric relation that we found in the $D_1$ case.
\begin{align*}
&(t^3x_1-x_4)(t^3x_1-x_5)(t^3x_1-x_6)\\ 
\equiv\,&(t^3x_1-x_4)(t^3x_1-x_5)(t^3x_1-x_6)\\
+\,& t^3x_1\left(t^6\sigma_2(x_1,x_2,x_3)-\sigma_2(x_4,x_5,x_6)\right)\\
-\, & t^6x_1^2\left(t^3(x_1+x_2+x_3)-(x_4+x_5+x_6)\right)\\
=\,&t^9x_1x_2x_3-x_4x_5x_6
\end{align*}
Let $\cP(D_2)=(t^9x_1x_2x_3-x_4x_5x_6)$ in $\cT_2^\prime$. Although technically defined as ideals in $\cT_i^\prime$, the $\cP(D_i)$ are contained in $\cR[x_0,\ldots,x_6,x_{13},\ldots,x_n]$, so we will abuse notation slightly in referring to them as $\cP=\cP(D_i)\subset\cT$.

For the ideals coming from $\cE^1(D_2)$, we have so far that the typical generator can be rewritten as $f_2\equiv(t^3x_1-x_4)w_{\text{in}}[x_5]+tx_9(t^{6+\weight}w_{\text{out}}[x_2]-w_{\text{in}}[x_5])$. Factoring out $t^2x_1-x_9$ yields
\begin{equation*}
\label{f2calc}
f_2\equiv (x_9-t^2x_1)(t^{6+\weight}w_{\text{out}}[x_2]-w_{\text{in}}[x_5])+(t^{8+\weight}w_{\text{out}}x_1[x_2]-t^{-1}w_{\text{in}}x_4[x_5]).
\end{equation*}
If $f_2$ is associated to a coherent region $R$ (assumed to contain elementary region $E_1$), then the second term above is an element of $\cE^{12}(D_2)$ associated to $V\cup E_2$, after eliminating $x_{11}$ and $x_{12}$ using the bivalent vertices on the left strand of $D_2$. Therefore, as an element of $\cT_2^\prime$, 
\begin{equation}
\label{finalnonlocal2}
f_2\equiv (x_9-t^2x_1)(t^{6+\weight}w_{\text{out}}[x_2]-w_{\text{in}}[x_5]).
\end{equation} 
This means that we should use polynomials of the form $t^{6+\weight}w_{\text{out}}[x_2]-w_{\text{in}}[x_5]$ as the generators for $\cE^x(D_2)$.

Using the expression for $f_2$ obtained in Equation~\ref{finalnonlocal2} and a computation of $x_9(x_9-t^2x_1)$ analogous to that preceding Equation~\ref{square} in the analysis of $x_9f_1$, we may express $x_9f_2$ as
\begin{align*}
x_9f_2
\equiv&\, (t^3x_1-x_5)(t^3x_1-x_6)(t^{6+\weight}w_{\text{out}}[x_2]-w_{\text{in}}[x_5]).
\end{align*}
This expression no longer contains any of $x_7,\ldots,x_{12}$. As in the $D_1$ case, we may use relations from $\cL^{123}_{\text{sym}}$ to rewrite this expression in a more convenient form. We exhibit the computation for the case when $x_2$ and $x_5$ do not appear. A similar method works for $t^{6+\weight}w_{\text{out}}x_2-w_{\text{in}}x_5$.
\begin{align*}
&\,(t^3x_1-x_5)(t^3x_1-x_6)(t^{6+\weight}w_{\text{out}}-w_{\text{in}})\\
\equiv&\,(t^3x_1-x_5)(t^3x_1-x_6)(t^{6+\weight}w_{\text{out}}-w_{\text{in}})\\
&+\,t^3x_1w_{\text{in}}(t^3x_1+t^3x_2+t^3x_3-x_4-x_5-x_6)\\
&-\,w_{\text{in}}(t^6x_1x_2+t^6x_1x_3+t^6x_2x_3-x_4x_5-x_4x_6-x_5x_6)\\
\equiv&\, (t^3x_1-x_5-x_6)(t^{9+\weight}w_{\text{out}}x_1-w_{\text{in}}x_4)+t^6(t^{\weight}w_{\text{out}}x_5x_6-w_{\text{in}}x_2x_3)
\end{align*}

The final expression above is a linear combination of relations that already hold in $\cT_2^\prime$. Suppose that $f_2$ comes from a coherent region $R$, which we assume contains $E_1$ and the elementary regions immediately above and below $E_2$, with associated subset $V_R$. Then the first relation, $t^{9+\weight}w_{\text{out}}x_1-w_{\text{in}}x_4$, is obtained by taking the union of $V_R$ with the 4-valent vertex in layer $s_2$ and the bivalent vertices immediately above and below it. Therefore, it is contained in $\cE^{12}(D_2)$. The second relation above, $t^{\weight}w_{\text{out}}x_5x_6-w_{\text{in}}x_2x_3$,  is associated to the subset obtained from $V_R$ by removing the 4-valent vertices in layers $s_1$ and $s_3$ and the bivalent vertex between edges $x_8$ and $x_{10}$. Therefore, it is contained in $\cN^\prime$.  As in the $D_1$ case, $x_9f_2$ contributes no new generators to $\cE(D_2)$ or $\cE^x(D_2)$.

Having computed the appropriate generators for $\cP(D_2)$, $\cP^x(D_2)$, $\cE(D_2)$, and $\cE^x(D_2)$, we have now established the following splitting of $\cA_{111}(D_2)$ as a direct sum of $\cR[x_0,\ldots,x_6,x_{13},\ldots,x_n]$-modules:
\[\cA_{111}(D_2)\isom\frac{\cT^\prime_2(1)}{\cP}\bigoplus\frac{\cT^\prime_2(t^2x_1-x_9)}{\cP^x(D_2)+\cE^x(D_2)}\]
Define the first summand to be $\cB(D_2)$ and the second summand to be $\cB^x(D_2)$.

\subsection*{Analysis of Summands and Edge Maps}
We claimed that the $\cB(D_i)$ are isomorphic as $\cR[x_0,\ldots,x_6,x_{13},\ldots,x_n]$-modules to $\cB$, the module defined at the beginning of Section~\ref{sec:reid3} and assigned to the diagram $D_\ast$ obtained from $D_i$ by replacing the region near the Reidemeister III move with the 6-valent vertex in Figure~\ref{6valentvertex}. 

Recalling that 
\[\cT=\frac{\mathcal{R}[x_0,\ldots,x_6,x_{13},\ldots,x_n]}{\cL^\prime+\cL^{123}_{\textrm{sym}}},\]
and comparing to the definition of $\cB$, we have a presentation of $\cB$ as
\[\cB\isom\frac{\cT}{(t^9x_1x_2x_3-x_4x_5x_6)+\cN^\prime+\cN^{123}_{\text{sym}}}.\]
Unwrapping the definitions of $\cT^\prime_i$, our final presentations of the $\cB(D_i)$ were
\begin{align*}
\cB(D_1)&\isom\frac{\cT}{\cN^\prime+\cE^{12}(D_1)+\cE^{13}(D_1)+\cE^{1234}(D_1)+\cP}\quad\quad\text{ and }\\
\cB(D_2)&\isom\frac{\cT}{\cN^\prime+\cE^{12}(D_2)+\cE^{123}(D_2)+\cE^{124}(D_2)+\cE^{1234}(D_2)+\cP}.
\end{align*}
Since $\cP=(t^9x_1x_2x_3-x_4x_5x_6)$, we only need to check that 
\begin{align*}
\cN^{123}_{\text{sym}}&=\cE^{12}(D_1)+\cE^{13}(D_1)+\cE^{1234}(D_1)\\
&=\cE^{12}(D_2)+\cE^{123}(D_2)+\cE^{124}(D_2)+\cE^{1234}(D_2).
\end{align*} This is done by comparing generating sets in Figure~\ref{6valnonloccompare}.

\begin{sidewaysfigure}
\begin{center}
\input{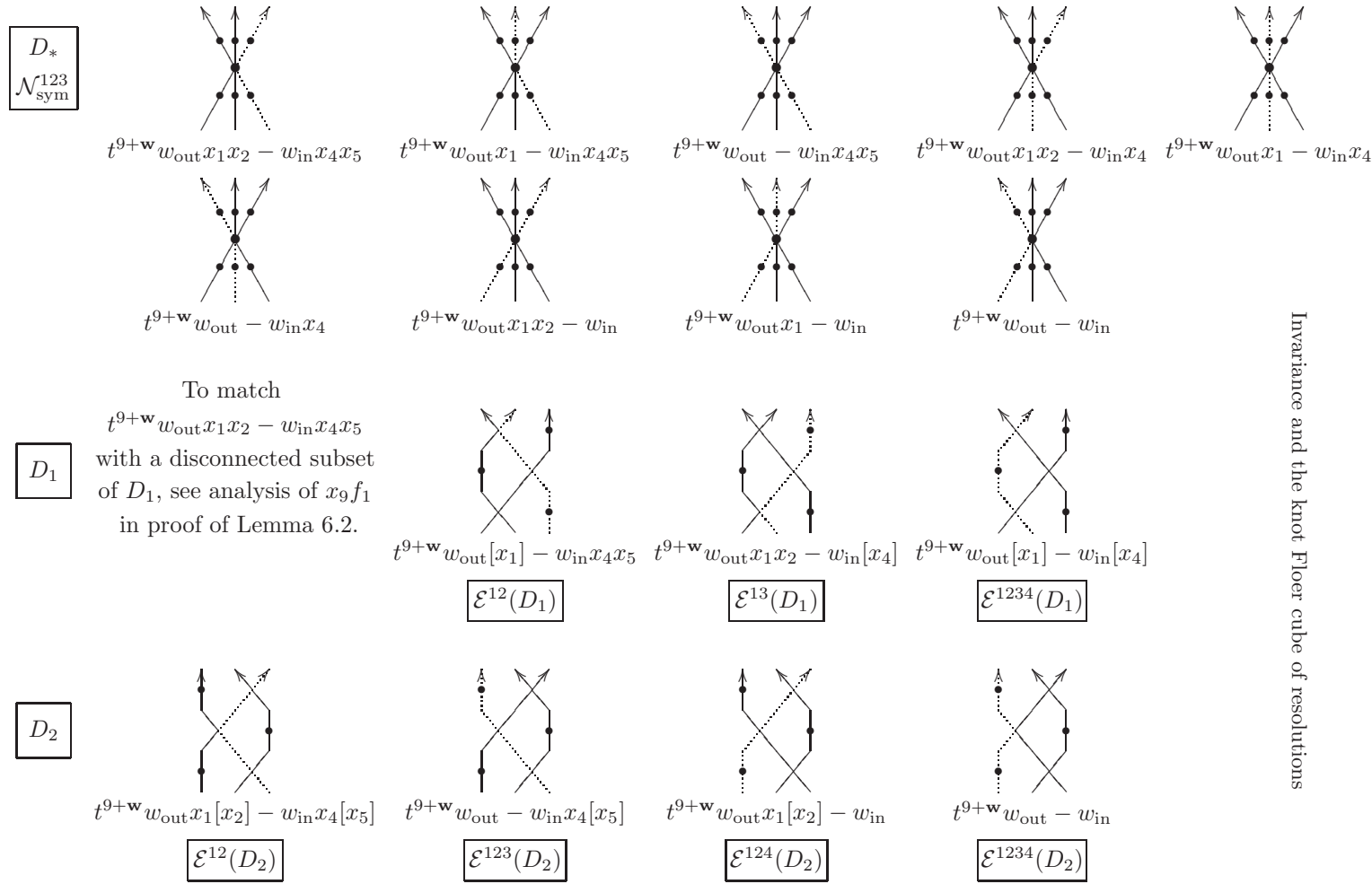}
\caption{The non-local relations generating $\cN^{123}_{\text{sym}}$ (shown in the first two rows) are the same as those associated to certain regions in the 111-resolutions of $D_1$ (row 3) and $D_2$ (row 4).  Assume that the braid axis is to the right of each picture. Brackets around $x_i$ denote an edge variable that may or may not occur in a relation depending on whether the coherent region under consideration contains an elementary region immediately above or below one of the elementary regions in the vicinity of the Reidemeister III move. Dotted lines show the boundary of the relevant coherent region when such adjacent elementary regions are not included. In each diagram, $\weight $, $w_{\text{out}}$, and $w_{\text{in}}$ come from the portions of the region not shown in these local pictures.}
\label{6valnonloccompare}
\end{center}
\end{sidewaysfigure}

It remains only to check that $\cB^x\isom\cC$ via the appropriate edge map. Our final presentations of the $\cB^x(D_i)$ were
\begin{align*}
\cB^x(D_1)&\isom\frac{\cT_1^\prime(t^2x_3-x_9)}{\cP^x(D_1)+\cE^x(D_1)}\quad\text{ and }\\
\cB^x(D_2)&\isom\frac{\cT_2^\prime(t^2x_1-x_9)}{\cP^x(D_2)+\cE^x(D_2)},
\end{align*}
while our final presentations of the $\cC(D_i)$ were
\begin{align*}
\cC(D_i)&\isom\frac{\cT(1)}{(r_i)+\cN_{101}(D_i)}
\end{align*}
Comparing to the notation in the proof of Lemma~\ref{reid3splitting101}, we see that $\cP^x(D_i)=(r_i)$, so the work is entirely in matching the generators of $\cN_{101}(D_i)$ with $\cE^x(D_i)$ and the various other ideals of non-local relations hidden in the definitions of the $\cT_i^\prime$. Specifically, we need to show that
\begin{align*}
\cN_{101}(D_1)&=\cE^x(D_1)+\cN^\prime+\cE^{12}(D_1)+\cE^{13}(D_1)+\cE^{1234}(D_1)\quad \text{ and }\\
\cN_{101}(D_2)&=\cE^x(D_2)+\cN^\prime+\cE^{12}(D_2)+\cE^{123}(D_2)+\cE^{124}(D_2)+\cE^{1234}(D_2).
\end{align*}

Since the 101- and 111-resolutions of $D_i$ differ only by whether the crossing in layer $s_2$ is singular or smooth, Observation~\ref{singsmoothobs} implies that $\cN_{111}(D_i)\subset\cN_{101}(D_i)$. Excepting $\cE^x(D_i)$, the ideals involved in the sums on the right hand side above are all contained in $\cN_{111}(D_i)$, so we have
\begin{align*}
\cN_{101}(D_1)&\supset\cN^\prime+\cE^{12}(D_1)+\cE^{13}(D_1)+\cE^{1234}(D_1)\quad \text{ and }\\
\cN_{101}(D_2)&\supset\cN^\prime+\cE^{12}(D_2)+\cE^{123}(D_2)+\cE^{124}(D_2)+\cE^{1234}(D_2).
\end{align*}
We can check directly that $\cN_{101}(D_i)\supset\cE^x(D_i)$ as well. Label the elementary regions in the 101-resolution of $D_i$ as in Figure~\ref{elemregions101}. Coherent regions in the 101-resolution of $D_1$ that use $F_1$ but not $F_2$ or $F_3$ correspond to coherent regions in the 111-resolution of $D_1$ that use $E_1$ but not any of the other $E_i$. Similarly, coherent regions in the 101-resolution of $D_2$ that use $F_1$ and $F_2$ but not $F_3$ correspond to coherent regions in the 111-resolution of $D_2$ that use $E_1$ but not any of the other $E_i$. In $\cN_{101}(D_i)$, such regions have associated non-local relations of the form $t^{3+\weight}w_{\text{out}}-w_{\text{in}}$ for $i=1$ or $t^{6+\weight}w_{\text{out}}[x_2]-w_{\text{in}}[x_5]$ for $i=2$ with the presence of $x_2$ and/or $x_5$ depending on whether the elementary regions immediately above and below $F_2$ are part of the coherent region under consideration. The corresponding regions in the 111-resolution have associated non-local relations in $\cE^x(D_i)$ which are exactly the same. Therefore, we have the full inclusions
\begin{align*}
\cN_{101}(D_1)&\supset\cE^x(D_1)+\cN^\prime+\cE^{12}(D_1)+\cE^{13}(D_1)+\cE^{1234}(D_1) \quad\text{ and }\\
\cN_{101}(D_2)&\supset\cE^x(D_2)+\cN^\prime+\cE^{12}(D_2)+\cE^{123}(D_2)+\cE^{124}(D_2)+\cE^{1234}(D_2).
\end{align*}

\begin{figure}[htbp]
\begin{center}
\input{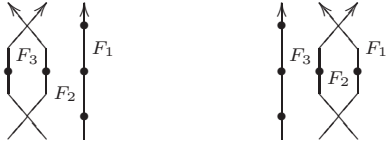}
\caption{Elementary regions in the 101-resolution of $D_1$ (left) and $D_2$ (right).}
\label{elemregions101}
\end{center}
\end{figure}

To prove the opposite inclusion, we use the coherent regions definition to classify the generators of $\cN_{101}(D_i)$ based on which of the $F_i$ are used, just as we did to understand $\cN_{111}$ earlier. Figure~\ref{111v101nonlocals} shows the general form of relations using various combinations of the $F_i$ in the 101-resolution of each $D_i$, and shows corresponding coherent regions in the 111-resolutions that produce the same non-local relations in one of the summands on the right-hand side of our desired equalities. 

\begin{sidewaysfigure}[htbp]
\begin{center}
\input{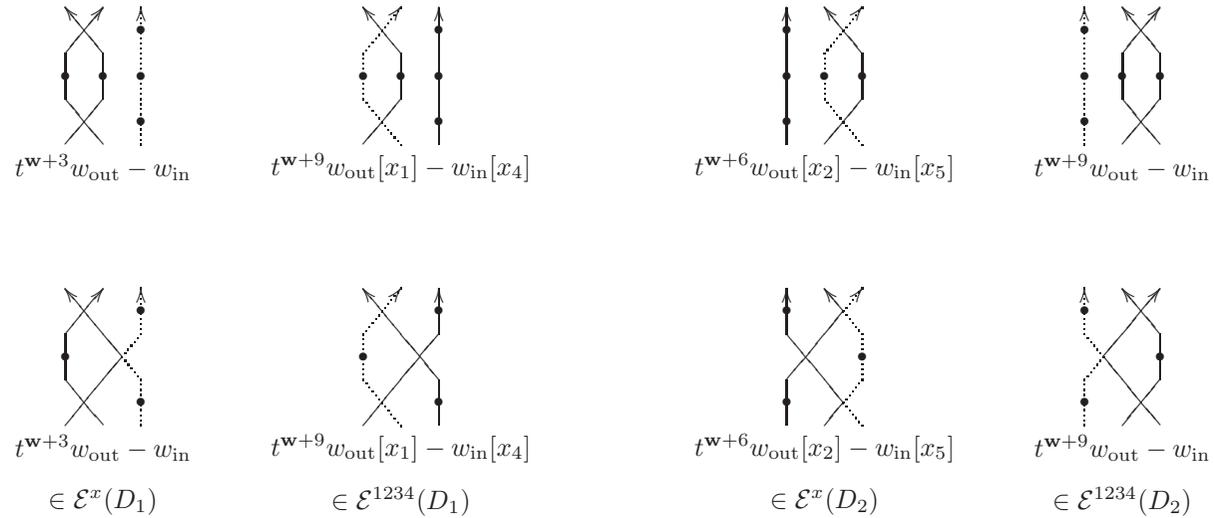}
\caption{The top row shows the general form of relations in $\cN_{101}(D_i)$ associated to coherent regions containing various combinations of the elementary regions $F_i$, with $\weight $, $w_{\text{out}}$, and $w_{\text{in}}$ coming from the portions of the region not shown in these local pictures. The bottom row identifies corresponding coherent regions in the 111-resolution of $D_i$ that have the same associated non-local relation.  Assume that the braid axis is to the right of each picture. Brackets around $x_i$ denote an edge variable that may or may not occur in a relation depending on whether the coherent region under consideration contains an elementary region immediately above or below one of the elementary regions in the vicinity of the Reidemeister III move. Dotted lines show the boundary of the relevant coherent region when such adjacent elementary regions are not included. }
\label{111v101nonlocals}
\end{center}
\end{sidewaysfigure}

We have now shown that $\cC$ and $\cB^x$ are quotients of $\cT$ by the same ideals. The edge map $\cA_{101}\to\cA_{111}$ is multiplication by $t^2x_3-x_9$ for $D_1$ and multiplication by $t^2x_1-x_9$ for $D_2$, so it maps the generator ($1$) of $\cA_{101}$ to the generator of $\cB^x$ in either case. This completes the proof that the edge map is an isomorphism when restricted to $\cC\to\cB^x$.
\end{proof}

\subsection{Simplifying Complexes for Reidemeister Move III}
\label{sec:reid3complexmanipulation}

Figures~\ref{reid3presentationsd1} and~\ref{reid3presentationsd2} show preferred generating sets for the ideals of local relations $\cL^{123}_{jk\ell}(D_i)$ associated to layers $s_1$, $s_2$, and $s_3$ in each resolution and preferred expressions for the differentials in $C(D_1)$ and $C(D_2)$. 

\begin{sidewaysfigure}
\begin{center}
\input{newreid3presentationsd1}
\caption{$C(D_1)$ as a complex of modules over $\cT[x_7,x_9]/\cN^\prime$ with preferred differentials. The column vectors below each diagram should be read as preferred generating sets for the ideal of local relations associated to the layers shown in that diagram. Choices made in constructing these preferred generating sets are explained at the beginning of Section~\ref{sec:reid3complexmanipulation}.}
\label{reid3presentationsd1}
\end{center}
\end{sidewaysfigure}

\begin{sidewaysfigure}
\begin{center}
\input{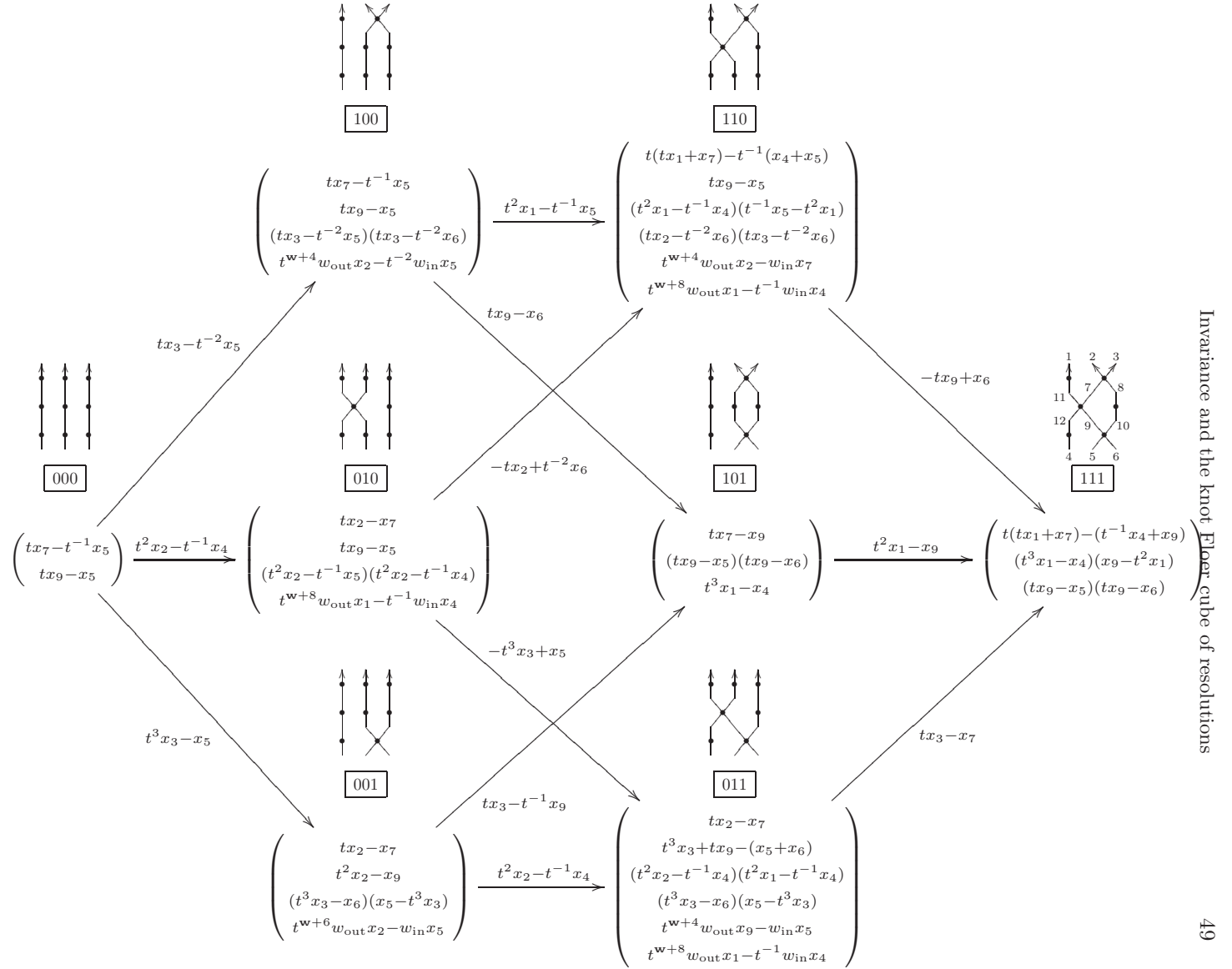}
\caption{$C(D_2)$ as a complex of modules over $\cT[x_7,x_9]/\cN^\prime$ with preferred differentials. The column vectors below each diagram should be read as preferred generating sets for the ideal of local relations associated to the layers shown in that diagram. Choices made in constructing these preferred generating sets are explained at the beginning of Section~\ref{sec:reid3complexmanipulation}.}
\label{reid3presentationsd2}
\end{center}
\end{sidewaysfigure}

The preferred generating sets should be viewed as presentations for the various $\cA_{jk\ell}(D_i)$ as modules over $\frac{\cT[x_7,x_9]}{\cN^\prime}$, where $\cN^\prime$ is generated by non-local relations associated to coherent regions not containing any of the elementary regions in the vicinity of the Reidemeister $\mathrm{III}$ move. We defined $\cT$ in Section~\ref{sec:reid3splittings} as $$\cT=\frac{\mathcal{R}[x_0,\ldots,x_6,x_{13},\ldots,x_n]}{\cL^\prime+\cL^{123}_{\textrm{sym}}},$$ where $\cL^\prime$ is generated by local relations associated to vertices away from the vicinity of the Reidemeister $\mathrm{III}$ move and $\cL^{123}_{\textrm{sym}}$ is generated by differences of certain elementary symmetric polynomials. In the process of proving the splitting lemmas of Section~\ref{sec:reid3splittings}, we established that the relations in $\cL^{123}_{\textrm{sym}}$ hold in $\cA_{101}(D_i)$ and $\cA_{111}(D_i)$. One can check that they hold for the other resolutions by comparatively straightforward manipulations of local relations. Since the relations in $\cL^\prime$, $\cL^{123}_{\text{sym}}$, and $\cN^\prime$ hold in all resolutions, we have omitted them from Figures~\ref{reid3presentationsd1} and~\ref{reid3presentationsd2}.

In building the preferred generating sets, we have eliminated $x_{10}$, $x_{11}$, and $x_{12}$ using the relations associated to bivalent vertices that appear in all of the resolutions: $tx_8-x_{10}$ in both $C(D_1)$ and $C(D_2)$; $tx_3-x_{11}$ and $tx_{12}-x_6$ in $C(D_1)$; and $tx_1-x_{11}$ and $tx_{12}-x_4$ in $C(D_2)$. Except in the 101- and 111-resolutions, it is possible to eliminate $x_7$, $x_8$, and $x_9$ as well in terms of $x_1,\ldots,x_6$. In these cases, we have listed linear relations used to do so for $x_7$ and $x_9$. Those used to eliminate $x_8$ can be inferred.

For all except the 101- and 111-resolutions, we have then listed quadratic relations in formats chosen to demonstrate the well-definedness of the isomorphism exhibited in Figure~\ref{reid3isom} and of incoming edge maps. Finally, we have listed representative non-local relations associated to coherent regions that use the elementary regions in the vicinity of the Reidemeister III move. 

For the 101- and 111-resolutions, we have listed the linear relations that were used to eliminate $x_7$ in the proofs of Lemmas~\ref{reid3splitting101} and~\ref{reid3splitting111}. We have then listed the quadratic relations that appear in the final presentations of these modules before the splitting step in the proofs of Lemmas~\ref{reid3splitting101} and~\ref{reid3splitting111}. We have omitted the non-local relations in these resolutions because they were already analyzed extensively in Section~\ref{sec:reid3splittings}.

The top lines of Figures~\ref{reid3startingcomplexesd1} and~\ref{reid3startingcomplexesd2} show condensed versions of $C(D_1)$ and $C(D_2)$ as presented in Figures~\ref{reid3presentationsd1} and~\ref{reid3presentationsd2}. The next step is to incorporate the splittings established in Section~\ref{sec:reid3splittings} and adjust matrix entries in the differentials accordingly. The bottom lines of Figures~\ref{reid3startingcomplexesd1} and~\ref{reid3startingcomplexesd2} show $C(D_1)$ and $C(D_2)$, respectively, after the splittings have been introduced and matrix entries adjusted. 

\begin{sidewaysfigure}
\begin{center}
\input{reid3startingcomplexesd1}
\caption{The top line is the complex $C(D_1)$, condensed from Figure~\ref{reid3presentationsd1} to emphasize differentials rather than module presentations. The bottom line incorporates the splittings from Section~\ref{sec:reid3splittings} into $C(D_1)$ and adjusts matrix entries in the differentials accordingly. Recall that $\cC^x(D_1)$ is generated by $tx_9-x_4$ and $\cB^x(D_1)$ is generated by $t^2x_3-x_9$.}
\label{reid3startingcomplexesd1}
\end{center}
\end{sidewaysfigure}

\begin{sidewaysfigure}
\begin{center}
\input{reid3startingcomplexesd2}
\caption{The top line is the complex $C(D_2)$, condensed from Figure~\ref{reid3presentationsd2} to emphasize differentials rather than module presentations. The bottom line incorporates the splittings from Section~\ref{sec:reid3splittings} into $C(D_2)$ and adjusts matrix entries in the differentials accordingly. Recall that $\cC^x(D_2)$ is generated by $tx_9-x_6$ and $\cB^x(D_2)$ is generated by $t^2x_1-x_9$.}
\label{reid3startingcomplexesd2}
\end{center}
\end{sidewaysfigure}

For $C(D_1)$ we adjust matrix entries as follows. In the rightmost matrix, we replace $x_7$ using the relation $t(tx_3+x_7)-(t^{-1}x_6-x_9)$, as we did when establishing the splitting of $\cA_{111}(D_1)$. We then arrange the entries in the rightmost matrix such that the row operation $\mathrm{I}_r+(t^2x_3-x_9)\mathrm{II}_r$ applied to the new matrix recovers the previous matrix. In the third column, this row operation produces $(t^2x_3-x_9)(tx_9-x_4)$, which is the image of the generator of $\cC^x$ in $\cA_{111}$. To see the equality directly, one must replace $x_9^2$ using the relation $(tx_9-x_5)(tx_9-x_4)$, as we did when establishing the splitting of $\cA_{111}(D_1)$ in Section~\ref{sec:reid3splittings}. For the middle matrix, we have arranged the entries such that applying the row operation $\mathrm{II}_r+(tx_9-x_4)\mathrm{III}_r$ to the new matrix recovers the previous matrix.

For $C(D_2)$, we handle the splitting in the same way. In the rightmost matrix, we replace $x_7$ in $\cA_{111}(D_2)$ using $t(tx_1+x_7)-(t^{-1}x_4+x_9)$ as we did when splitting $\cA_{111}(D_2)$. Then, we arrange matrix entries such that applying the row operation $\mathrm{I}_r+(t^2x_1-x_9)\mathrm{II}_r$ to the new matrix on the right recovers the previous matrix on the right. This entails replacing $x_9^2$ using the relation $(tx_9-x_5)(tx_9-x_6)$ as we did when establishing the splitting of $\cA_{111}(D_2)$. In the middle matrix, the row operation $\mathrm{II}_r+(tx_9-x_6)\mathrm{III}_r$ applied to the new matrix recovers the previous matrix.

We are now prepared to perform the changes of basis necessary to identify contractible summands in $C(D_1)$ and $C(D_2)$. Figure~\ref{reid3changebasis} exhibits the operations step by step in the case of $C(D_1)$. We may use almost the same row and column operations in the case of $C(D_2)$. Simply exchange $x_4$ with $x_6$ and $x_1$ with $x_3$, but leave everything else the same. For example, the second change of basis performed on $C(D_1)$ used row operation $\mathrm{I}_r-(t^2x_3-t^{-1}x_5)\mathrm{III}_r$ on the middle matrix. The second change of basis on $C(D_2)$ should instead use row operation $\mathrm{I}_r-(t^2x_1-t^{-1}x_5)\mathrm{III}_r$ on the middle matrix.

\begin{sidewaysfigure}
\begin{center}
\input{reid3changebasis}
\caption{Beginning with $C(D_1)$ as presented in the bottom line of Figure~\ref{reid3startingcomplexesd1}, we change basis several times. The change of basis occurs in the boxed homological grading, with corresponding row and column operations on incoming and outgoing maps indicated below the arrows. The result is a complex with contractible summands $\cA_{100}\to\cC^x$ and $\cC\to\cB^x$.}
\label{reid3changebasis}
\end{center}
\end{sidewaysfigure}

Figures~\ref{reid3removesummandsd1} and~\ref{reid3removesummandsd2} show $C(D_1)$ and $C(D_2)$ after the changes of basis are complete (top two lines), and then the complexes obtained by removing the contractible summands $\cA_{100}(D_i)\to\cC^x(D_i)$ and $\cC(D_i)\to\cB^x(D_i)$ (bottom two lines). The bottom two complexes are $\overline{C}(D_1)$ and $\overline{C}(D_2)$ as described at the beginning of Section~\ref{sec:reid3}. We have now indicated the induced edge maps as well.

\begin{sidewaysfigure}
\begin{center}
\input{reid3removesummandsd1}
\caption{The top line shows $C(D_1)$ after the changes of basis shown in Figure~\ref{reid3changebasis}. The bottom line shows the result of removing contractible summands $\cA_{100}\to\cC^x$ and $\cC\to\cB^x$, which is $\overline{C}(D_1)$ as described at the beginning of Section~\ref{sec:reid3}.}
\label{reid3removesummandsd1}
\end{center}
\end{sidewaysfigure}

\begin{sidewaysfigure}
\begin{center}
\input{reid3removesummandsd2}
\caption{The top line shows $C(D_2)$ after changes of basis analogous to those shown in Figure~\ref{reid3changebasis} for $\overline{C}(D_1)$. The bottom line shows the result of removing contractible summands $\cA_{100}\to\cC^x$ and $\cC\to\cB^x$, which is $\overline{C}(D_2)$ as  described at the beginning of Section~\ref{sec:reid3}.}
\label{reid3removesummandsd2}
\end{center}
\end{sidewaysfigure}

The last step is simply to exhibit an isomorphism between the simplified complexes we have obtained for $C(D_1)$ and $C(D_2)$. Figure~\ref{reid3isom} shows the appropriate chain map from $\overline{C}(D_1)$ to $\overline{C}(D_2)$. The map identifies the summand $\cA_{010}(D_1)$ with $\cA_{001}(D_2)$; $\cA_{001}(D_1)$ with $\cA_{010}(D_2)$; $\cA_{110}(D_1)$ with $\cA_{011}(D_2)$; and $\cA_{011}(D_1)$ with $\cA_{110}(D_2)$. It also multiplies some summands by a power of $t$, which (1) directly maps the relevant presentations from Figure~\ref{reid3presentationsd1} to the corresponding presentations from Figure~\ref{reid3presentationsd2}; (2) accounts for the differences in exponents of $t$ in the edge maps of the simplified complexes; and (3) effectively moves layers of marked points past layers of 4-valent vertices to match the resolutions of $D_1$ with the resolutions of $D_2$ that remain in the simplified complexes. Finding the inverse of the map in Figure~\ref{reid3isom} is a straightforward computation.

\begin{sidewaysfigure}
\begin{center}
\input{reid3isom}
\caption{An isomorphism from $\overline{C}(D_1)$ (top) to $\overline{C}(D_2)$ (bottom).}
\label{reid3isom}
\end{center}
\end{sidewaysfigure}

\section{Conjugation}
\label{sec:movestar}

\begin{figure}[tbp]
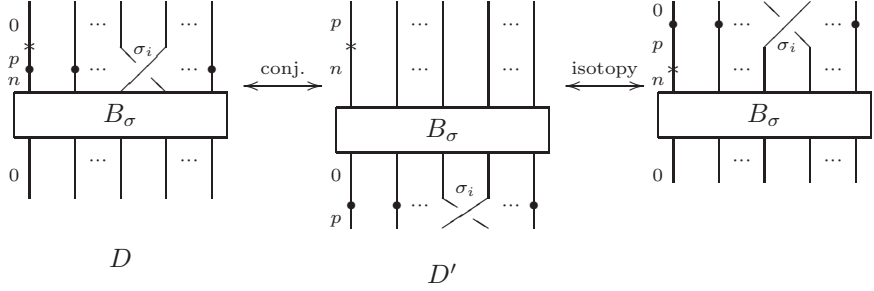

\begin{center}
$$\xymatrix@+5pt{{\scalebox{1}{\input{starbeforebiv}}}\ar@{<->}[r]^{\text{conj.}}
& {\scalebox{1}{\input{starafterbiv}}}\ar@{<->}[r]^{\text{isotopy}}
& {\scalebox{1}{\input{starbiviso}}}
}$$
\caption{Diagram for Lemma~\ref{lemma:movestarbiv}: moving the basepoint across a bivalent vertex is equivalent to conjugating $\sigma_i\sigma$ to $\sigma\sigma_i$ for $i\neq b-1$. These diagrams should be viewed as closed braids, although the closure strands are not shown. Assume that all strands are oriented upwards and that the braid axis is to the right of each picture.}
\label{bivstarmove}
\end{center}
\end{figure}

\begin{figure}[tbp]
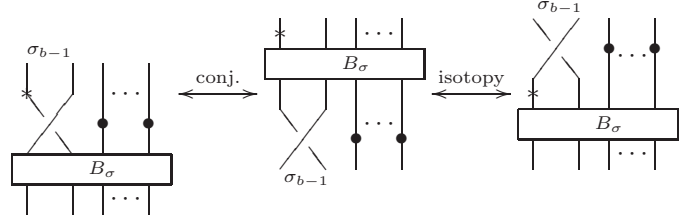

	\begin{center}
		$$\xymatrix@+5pt{\input{starcrossing1.tex} \ar@{<->}[r]^{\text{conj.}} & 
\input{starcrossing6.tex} \ar@{<->}[r]^{\text{isotopy}} &
\input{starcrossing5.tex}}
		$$
	\caption{Conjugating $\sigma_{b-1}\sigma$ to $\sigma\sigma_{b-1}$ would be equivalent to moving the basepoint across $\sigma_{b-1}$. These diagrams should be viewed as closed braids, although the closure strands are not shown. Assume that all strands are oriented upwards and that the braid axis is to the right of each picture.}
		\label{movestarcrossing}
	\end{center}
\end{figure}

In this section we demonstrate that $\mathcal{A}_I(D)$ is invariant under conjugation of the layered braid diagram $D$ away from the basepoint. That is, conjugation is allowed only by braid generators $\sigma_1,\ldots,\sigma_{b-2}$ and not by $\sigma_{b-1}$. (Our convention is to label the braid generators from right to left.) This limitation arises because the basepoint has a role in determining which cycles, subsets, or regions are used to define non-local relations. (Since conjugation is a planar isotopy of a braid diagram, it does not change the edge ring or the local relations.)  Proving that $\cA_I(D)$ is invariant under conjugation by one of $\sigma_1,\ldots,\sigma_{b-2}$ is equivalent to proving that it is invariant under moving the basepoint past a bivalent vertex, as in Figure~\ref{bivstarmove}. We do so in Lemma~\ref{lemma:movestarbiv}. Proving that $\cA_I(D)$ is invariant under conjugation by $\sigma_{b-1}$ would require moving the basepoint past a crossing, as in Figure~\ref{movestarcrossing}. We have so far been unable to do this.

\begin{lemma}
\label{lemma:movestarbiv}
Let $D$ be the layered braid diagram for a braid word of the form $\sigma_i\sigma$, where $i\neq b-1$ and $\sigma$ is any braid word. Let $D^\prime$ be the layered braid diagram for $\sigma\sigma_i$. Fix edge labels as in Figure~\ref{bivstarmove} with $p>n$. Then for any index $I$, $\cA_I(D)\isom\cA_I(D^\prime)$ as $\cR[x]$-algebras, where $x$ acts as the variable associated to the vertex outgoing from the basepoint in each diagram.
\end{lemma}

\begin{proof}
Whether $\sigma_i$ is resolved or singularized in the $I$-resolution of $D$, Figure~\ref{bivstarmove} indicates that it suffices to prove that we can move the basepoint across a bivalent vertex on the leftmost strand. Let $x_0,\ldots,x_n,x_p$ denote the variables in the edge ring for $D$ and $y_0,\ldots,y_n,y_p$ denote the variables in the edge ring for $D^\prime$. Edge labels are shown in Figure~\ref{bivstarmove} for edges $x_i$ and $y_i$ when $i\in\{0,n,p\}$. The remaining edges are labeled such that the position of $x_i$ in the diagram on the left of Figure~\ref{bivstarmove} matches the position of $y_i$ in the diagram on the right of Figure~\ref{bivstarmove}.

We will view $\cA_I(D)$ and $\cA_I(D^\prime)$ as $\cR[x]$-modules by equating $x$ with $x_0$ and with $y_p$, respectively. Define an $\cR[x]$-module map $\varphi:\cA_I(D)\to\cA_I(D^\prime)$ by $x_0\mapsto y_p$, $x_i\mapsto ty_i$ for $1\leq i\leq n$, and $x_p\mapsto y_n$. To see that it is well-defined and an isomorphism, first notice that $\varphi$ maps the linear relation $tx_p-x_n$ (coming from the bivalent vertex nearest the basepoint in $D$) to 0 in $\cA_I(D^\prime)$. Now use the relation $tx_p-x_n$ to find a presentation of $\cA_I(D)$ in which $x_p$ does not appear. Suppose that $f(x_0,\ldots,x_n)$ is one of the relations in this presentation. Since all of the relations in the original presentation of $\cA_I(D)$ were homogenous in the $x_i$, $f$ will be as well. Then $\varphi(f(x_0,\ldots,x_n))=f(y_p,ty_1,\ldots,ty_n)\equiv f(ty_0,ty_1,\ldots,ty_n)\equiv f(y_0,y_1,\ldots,y_n),$ where ``$\equiv$'' here means ``generates the same ideal in $\cR[y_0,\ldots,y_n,y_p]/(ty_0-y_p)$.'' Since $ty_0-y_p$ is a relation in $\cA_I(D^\prime)$ (associated to the bivalent vertex nearest the basepoint), this calculation says that $\varphi$ identifies each relation in the chosen presentation of $\cA_I(D)$ with a relation in $\cA_I(D^\prime)$. The map defined by $y_i\mapsto t^{-1}x_i$ for ($0\leq i\leq n$) and $y_p\mapsto x_0$ is an inverse for $\varphi$, which one can check is well-defined by a similar argument. 
\end{proof}

\section{Stabilization / Reidemeister Move I}
\label{sec:reid1}

Throughout this section, we work over $\widehat{\cR}[\underline{x}(D)]$. The arguments presented here require the completion of the ground ring because we invert an element of the form $1-t^k$, and because the completion is required to make Observation~\ref{disconnectedobs} about disconnected resolutions hold.

\begin{figure}[tbp]
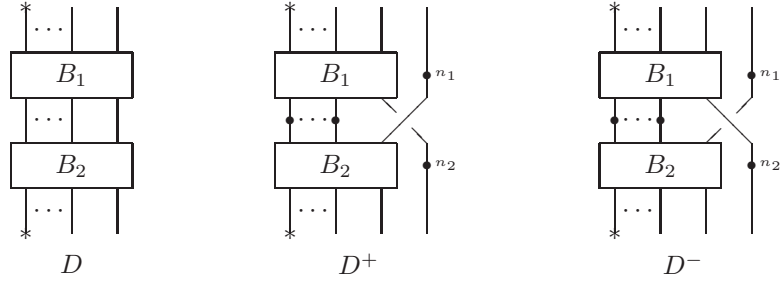

	\begin{center}
		\subfloat{\input{reid1origin.tex}}\quad\quad\quad\quad\quad
		\subfloat{\input{reid1plusin.tex}}\quad\quad\quad\quad\quad
		\subfloat{\input{reid1minusin.tex}}
		\caption{Diagrams $D$, $D^{+}$, and $D^{-}$.These diagrams should be viewed as closed braids, although the closure strands are not shown. Assume that all strands are oriented upwards and that the braid axis is to the right of each picture.}
	\label{innerreid1}
	\end{center}
\end{figure}

\begin{figure}[tbp]
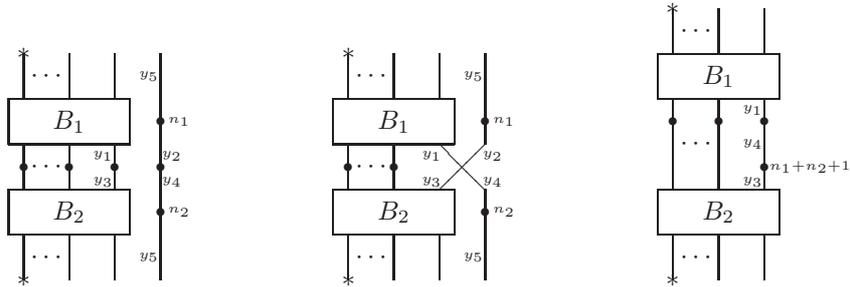

	\begin{center}
		\subfloat{\input{reid1smooth.tex}}\quad\quad\quad\quad\quad
		\subfloat{\input{reid1sing.tex}}\quad\quad\quad\quad\quad
		\subfloat{\input{reid1markedpt.tex}}
		\caption{From left to right: the smoothed resolution of layer $s$ in $D^{+}$ of $D^{-}$; the singular resolution of layer $s$ in $D^{+}$ of $D^{-}$; the diagram $D^{\bullet}$. These diagrams should be viewed as closed braids, although the closure strands are not shown. Assume that all strands are oriented upwards and that the braid axis is to the right of each picture.}
	\label{innerreid1resolutions}
	\end{center}
\end{figure}

Let $D$ and $D^{+}$ (respectively $D^{-}$) be closed braid projections that differ by a positive (respectively negative) stabilization in layer $s$ on the innermost strand as in Figure~\ref{innerreid1}. Ideally, we would like to show that $C(D)$, $C(D^{+})$, and $C(D^{-})$ are chain homotopy equivalent. Figure~\ref{innerreid1resolutions} shows the two possible resolutions of the crossing on the innermost strand. Any resolution in which the crossing on the innermost strand is smoothed is disconnected, so the corresponding algebra $\cA_{\mathrm{I}_11\mathrm{I}_2}(D^{+})$ or $\cA_{\mathrm{I}_10\mathrm{I}_2}(D^{-})$ vanishes by Observation~\ref{disconnectedobs}. It would suffice, then, to show that $\cA_{\mathrm{I}_10\mathrm{I}_2}(D^{+})\isom\cA_{\mathrm{I}_11\mathrm{I}_2}(D^{-})\isom\cA_{\mathrm{I}_1\mathrm{I}_2}(D)$ for any resolutions $\mathrm{I}_j\in\{0,1\}^{n_j}$.

However, the behavior of the algebras under stabilization is not so straightforward. Instead, the algebras associated to the resolutions in which the crossing on the innermost strand is singularized are isomorphic to the algebras associated to the corresponding resolutions of $D^{\bullet}$ (shown on the right in Figure~\ref{innerreid1resolutions}). In other words, stabilizing once on the innermost strand is equivalent to adding a marked point to the edge where the stabilization occurs. The marked point has weight equal to the total weight of the strand being added, plus one. Equivalently, it has weight equal to the total number of crossings in the diagram, plus one.

The diagram $D^\bullet$ is not a layered braid, but the definition of $\cA$ extends easily to encompass this case. Use the same local relations for crossings and other bivalent vertices, along with a relation $t^{n_1+n_2+1}y_4-y_3$ associated to the new marked point. Define the non-local relations as before, but adjust their weights upwards by $n_1+n_2+1$ if they are associated to a subset, cycle, or region encompassing the new marked point. 

With this definition, $\cA_I(D^\bullet)$ is still isomorphic to singular knot Floer homology with twisted coefficients, as described in Section~\ref{sec:algdfn}. It still fits into the skein exact sequence described in~\cite{ozsszcube}. (These facts are not proved here, but would follow from arguments similar to those in~\cite{ozsszcube} and those in Section~\ref{sec:computeshfk}.) Unfortunately, our proofs of the categorified MOY relations underlying Reidemeister moves II and III do not extend to diagrams like $D^\bullet$ that have bivalent vertices with different weights on different strands. It may be that a more subtle version of the MOY calculus, taking into account bivalent vertices of various weights, could unify our proofs of the categorified braid-like MOY moves with our description of stabilization. We are not currently aware of such a model.

\begin{prop}
\label{reid1prop}
Let $D^{+}$ and $D^{-}$ be the diagrams in Figure~\ref{innerreid1} with $n_i$ crossings in $B_i$. Figure~\ref{innerreid1resolutions} shows the $(\mathrm{I}_10\mathrm{I}_2)$-resolution of $D^{+}$, which is identical to the $(\mathrm{I}_11\mathrm{I}_2)$-resolution of $D^{-}$. Let $D^{\bullet}$ be the diagram on the right in Figure~\ref{innerreid1resolutions} with edge labels as shown. Let $x_0,\ldots,x_k$ be the edges in the unlabeled portion of all of these diagrams. Then there are isomorphisms of $\widehat{\cR}[x_0,\ldots,x_k,y_1,y_3,y_4]$-modules
\[\cA_{\mathrm{I}_10\mathrm{I}_2}(D^{+})\isom\cA_{\mathrm{I}_11\mathrm{I}_2}(D^{-})\isom\cA_{\mathrm{I}_1\mathrm{I}_2}(D^{\bullet})\] for all resolutions $\mathrm{I}_j\in\{0,1\}^{n_j}$.
\end{prop}

\begin{proof}
The $(\mathrm{I}_10\mathrm{I}_2)$-resolution of $D^{+}$ is identical to the $(\mathrm{I}_11\mathrm{I}_2)$-resolution of $D^{-}$, so we will refer to the $(\mathrm{I}_10\mathrm{I}_2)$-resolution of $D^{+}$ throughout this proof without loss of generality. The edge ring of the $(\mathrm{I}_10\mathrm{I}_2)$-resolution of $D^{+}$ is $\widehat{\cR}[x_0,\ldots,x_k,y_1,y_2,y_3,y_4,y_5]$.  Let $n=n_1+n_2$. The marked points on the innermost strand give relations $t^{n_2}y_4-y_5$ and $t^{n}y_4-y_2$. Use these relations to eliminate $y_2$ and $y_5$ from the edge ring, leaving $\widehat{\cR}[x_0,\ldots,x_k,y_1,y_3,y_4]$, which is the edge ring for $D^{\bullet}$. We will work in the context of $\widehat{\cR}[x_0,\ldots,x_k,y_1,y_3,y_4]$-modules for the remainder of this proof.

The local relations for crossings not on the innermost strand do not use $y_2$, $y_4$, or $y_5$, and are the same for the singular resolutions of $D^{+}$ and $D^{-}$ as they are for $D^{\bullet}$. Let $\cL$ be the ideal generated by these local relations.

Non-local relations associated to coherent regions in $D^{\bullet}$ have the form \[t^{w^\prime+n+2}w_{\text{out}}-w_{\text{in}},\] where $y_1$, $y_3$, and $y_4$ do not divide $w_{\text{out}}$ or $w_{\text{in}}$, and $w^\prime$ is the contribution to the region's weight from vertices other than those pictured on the innermost strand.  Let $\cN$ be the ideal generated by such relations.  Each of these relations corresponds to a non-local relation in the singular resolution of $D^{+}$ and $D^{-}$ as follows. Let $E_1$, $E_2,\ldots,E_p$ be the elementary regions in $D^{\bullet}$ with $E_1$ the region containing the braid axis. Then the elementary regions in the singular resolution of $D^{+}$ and $D^{-}$ are $E_2,\ldots,E_p$ along with two others: the region containing the braid axis, which we will call $E_a$, and the region adjacent to $E_a$, which we will call $E_b$. Any coherent region in $D^{\bullet}$ contains $E_1$, so can be written as $R=E_1\cup E_{i_1}\cup \cdots\cup E_{i_r}$.  The region $E_a\cup E_b \cup E_{i_1}\cup \cdots\cup E_{i_r}$ in the singular resolution of $D^{+}$ and $D^{-}$ has the same associated non-local relation as $R$ does.

The relations defining $\cA_{\mathrm{I}_10\mathrm{I}_2}(D^{+})\isom
\cA_{\mathrm{I}_11\mathrm{I}_2}(D^{-})$ that we have not yet accounted for are as follows. (Recall that we have eliminated $y_2$ using the relation $t^{n}y_4-y_2$.)
\[\cJ=\left(\begin{array}{c}
ty_1+t^{n+1}y_4-y_3-y_4\\
t^{n+2}y_1y_4-y_3y_4\\
t^{n+2}y_1-y_3
\end{array}\right)\]
The first two lines come from local relations associated to the singular crossing on the innermost strand. The third line is the non-local relation associated to the coherent region $E_a$. So far, we have established that
\[\cA_{\mathrm{I}_10\mathrm{I}_2}(D^{+})\isom
\cA_{\mathrm{I}_11\mathrm{I}_2}(D^{-})\isom
\frac{\widehat{\cR}[x_0,\ldots,x_k,y_1,y_3,y_4]}{\cL+\cN+\cJ}
\]

We now simplify the presentation of $\cJ$. Perform the row operation $\mathrm{I}-\mathrm{III}$ to transform the first line into $(t^{n+1}-1)(y_4-ty_1)$. Since $t^{n+1}-1$ is a unit in $\widehat{\cR}$, we may remove that factor. Then factor out $y_4$ from the second line in the presentation of $\cJ$ to see that the second line is a multiple of the first, hence can be discarded. we then obtain the following alternative presentation of $\cJ$.
\[\cJ=\left(\begin{array}{c}
y_4-ty_1\\
t^{n+2}y_1-y_3
\end{array}\right)\]

The generators in this alternative presentation of $\cJ$ can also be viewed as the local relations associated to marked points visible on the innermost strand of $D^{\bullet}$, which are exactly the relations defining
$\cA_{\mathrm{I}_1\mathrm{I}_2}(D^{\bullet})$ that we had not yet accounted for above. That is, 
\[\cA_{\mathrm{I}_1\mathrm{I}_2}(D^{\bullet})\isom
\frac{\widehat{\cR}[x_0,\ldots,x_k,y_1,y_3,y_4]}{\cL+\cN+\cJ}.
\]
\end{proof}

\section{Identification with knot Floer homology}
\label{sec:computeshfk}

The set-up of the cube of resolutions in Section~\ref{sec:definitions} of this paper differs somewhat from
 Ozsv\'ath and Szab\'o's original formulation \cite{ozsszcube}, so it does not follow formally from their work that $C(D)$, as defined in (\ref{cubedfn}) of this paper, computes knot Floer homology. However, an adaptation of the arguments in Sections 3--5 of \cite{ozsszcube}, suffices to prove the following result, which is an analogue of \cite[Theorem~1.2]{ozsszcube}.

\begin{prop}
\label{prop:computeshfk}
Let $D$ be a layered braid diagram with initial edge $x_0$. Then there is an
isomorphism of graded $\mathbb{Z}_2[x_0]$-modules
$$H_\ast(C(D)\otimes_{\cR[\underline{x}(D)]}\widehat{\cR[\underline{x}(D)]}\otimes\mathbb{Z}_2)\isom \hfkmin(K)\otimes_{\mathbb{Z}_2}\mathbb{Z}_2[t^{-1},t]]$$
and an isomorphism of graded $\mathbb{Z}_2$-vector spaces
$$H_\ast(C(D)/(x_0)\otimes_{\cR[\underline{x}(D)]}\widehat{\cR[\underline{x}(D)]}\otimes\mathbb{Z}_2)\isom \hfkhat(K)\otimes_{\mathbb{Z}_2}\mathbb{Z}_2[t^{-1},t]].$$
\end{prop}

The two key differences between our set-up and that of \cite{ozsszcube} are the use of layered braid diagrams and the ground ring over which we define the cube of resolutions chain complex. Ozsv\'ath and Szab\'o use a knot projection in braid form with a basepoint $\ast$, but do not require the additional bivalent vertices that we add parallel to each crossing when forming a layered braid diagram.  Consequently, in their diagrams, bivalent vertices arise only when a crossing is smoothed, which means they come in pairs that lie on adjacent strands. A layered braid diagram has these sorts of bivalent vertices, but also others. This difference will require us to modify the Heegaard diagrams used in the proof of \cite[Theorem~1.2]{ozsszcube}.

The second difference between our set-up and that of \cite{ozsszcube} is in the ground rings over which the cube of resolutions complexes are defined. We define the algebras $\cA_I(D)$ over $\cR[\underline{x}(D)]=\mathbb{Z}[t^{-1},t][\underline{x}(D)]$, and pass to the completion $\widehat{\cR}[\underline{x}(D)]=\mathbb{Z}[t^{-1},t]][\underline{x}(D)]$ only when describing the behavior of the chain complex under stabilization. Ozsv\'ath and Szab\'o set up their algebras over $\mathbb{Z}_2[\underline{x}(D),t]$, then pass to $\mathbb{Z}_2[\underline{x}(D)][t^{-1},t]]$ to identify the homology of their cube of resolutions chain complex with knot Floer homology~\cite[Theorem~1.2]{ozsszcube}. Their algebras are defined as the singular knot Floer homology with particular choices of twisted coefficients. They require power series in $t$ with coefficients in $\mathbb{Z}_2[\underline{x}(D)]$ to make the singular knot Floer homology well-defined (by ensuring that its differential is a finite sum). They need to invert $t$ to apply their Lemma 2.2, which shows that knot Floer homology with twisted coefficients is isomorphic to the usual knot Floer homology tensored with an extended ground ring. These choices of rings in each case allow results to be stated in the greatest possible generality, but a profusion of tensor products will be required to bring the two approaches into alignment.

\begin{proof}
Ozsv\'ath and Szab\'o prove \cite[Theorem~1.2]{ozsszcube} in three steps: calculate a particular twisting of singular knot Floer homology to verify that it is identical to the algebra they define as a quotient of the edge ring \cite[Section~3]{ozsszcube}; establish a spectral sequence from the cube of resolutions defined algebraically to knot Floer homology \cite[Section 4]{ozsszcube}; show that the spectral sequence collapses \cite[Section 5]{ozsszcube}. We mirror each of these arguments in turn, pointing out where modifications are required to address the differences between our set-up (Section~\ref{sec:definitions} of this paper) and that of \cite{ozsszcube}.

Let $S$ be a layered braid diagram with all crossings singularized or smoothed. The twisted version of singular knot Floer homology needed to recover the algebra $\cA(S)$ as defined in~(\ref{corneralg}) in Section~\ref{sec:algdfn} of this paper is specified by the ``initial diagram'' in \cite[Figure~3]{ozsszcube} with the additional rule that near a bivalent vertex that does not arise from smoothing a crossing, the diagram has the form shown on the left in Figure~\ref{fig:hgddiagram}. Near a pair of bivalent vertices that arise from smoothing a crossing, we use the same diagram as in \cite[Figure~3]{ozsszcube}, which is shown in the middle in Figure~\ref{fig:hgddiagram}. Call this the \emph{modified initial diagram}. Let $\underline{CFK}^-(S)$ denote the chain complex coming from the modified initial diagram. That is, $\underline{CFK}^-(S)$ is the $\mathbb{Z}_2[\underline{x}(S)][[t]]$-module whose generators are given by intersection points and differentials by counting holomorphic disks with respect to the twisting in the modified initial diagram. See~\cite{ozsszstipsing} for a precise definition of singular knot Floer homology, \cite[Section~2.1]{ozsszcube} for details on twisted coefficients in knot Floer homology generally, and \cite[Section~3.1]{ozsszcube} for details on combining singular knot Floer homology with twisted coefficients. The completion of the ground ring with respect to $t$ is necessary to make the differential on twisted singular knot Floer homology well defined, as detailed in \cite[Section~3.1]{ozsszcube}. We will continue to work over $\mathbb{Z}_2[\underline{x}(S)][[t]]$ for the first section of this proof, so abbreviate this ring by $\cR^\prime$. 

\begin{figure}[tbp]
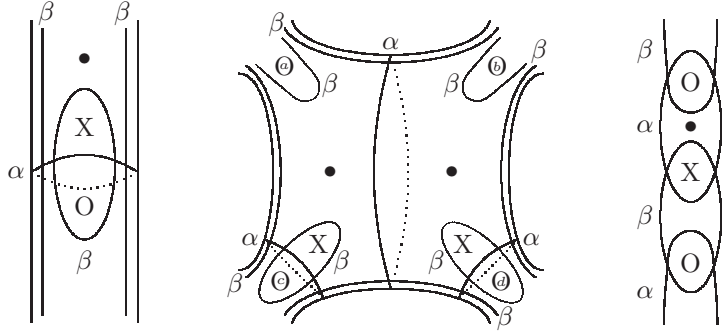

\begin{center}
\subfloat{\scalebox{1}{\input{initialextra.tex}}}\quad \quad \quad
\subfloat{\scalebox{1}{\input{initialsmoothed.tex}}}\quad \quad \quad
\subfloat{\scalebox{1}{\input{planar.tex}}}
\caption{From left to right: the modified initial diagram near an extra bivalent vertex; the modified initial diagram near a bivalent vertex arising from a smoothing; the planar diagram or the master diagram near any bivalent vertex. The bold dots in each picture show the marking that specifies our particular twisted version of singular knot Floer homology.}
\label{fig:hgddiagram}
\end{center}
\end{figure}

Let $M$ denote the Koszul complex on the linear relations for each vertex.
$$M=\bigotimes_{v\in V_4}\left(\cR^\prime\xrightarrow{tx_a^{(v)}+tx_b^{(v)}-x_c^{(v)}-x_d^{(v)}}\cR^\prime\right)\otimes\bigotimes_{v\in V_2}\left(\cR^\prime\xrightarrow{tx_a^{(v)}-x_c^{(v)}}\cR^\prime\right),$$
where $V_4$ and $V_2$ denote the set of 4-valent and bivalent vertices, respectively, in $S$. Let $C^\prime(S)=\underline{CFK}^-(S)\otimes M$. Then the claim, an analogue of \cite[Theorem~3.1]{ozsszcube}, is that we can identify $H_\ast(C^\prime(S))$ with $\cA(S)$ after appropriately changing the ground rings.  Recall that $\cA(S)$ was defined in~(\ref{corneralg}) of Section~\ref{sec:algdfn} of this paper as an $\cR[\underline{x}(S)]=\mathbb{Z}[t^{-1},t][\underline{x}(S)]$-module. Therefore, the precise claim is that
\begin{equation}
\label{singcalc}
H_\ast\!\left(C^\prime(S)\right)\otimes_{\cR^\prime}\cR^\prime[t^{-1}]\isom \cA(S)\otimes_{\cR[\underline{x}(S)]}\widehat{\cR[\underline{x}(S)]}\otimes\mathbb{Z}_2.
\end{equation}
The reduced version of the statement, 
\begin{equation}
\label{singcalcreduced}
H_\ast\!\left(C^\prime(S)/(x_0)\right)\otimes_{\cR^\prime}\cR^\prime[t^{-1}]\isom \cA(S)/(x_0)\otimes_{\cR[\underline{x}(S)]}\widehat{\cR[\underline{x}(S)]}\otimes\mathbb{Z}_2,
\end{equation}
then follows immediately.

The arguments required to prove \cite[Proposition~3.4]{ozsszcube} apply essentially unchanged to show that $H_\ast(C^\prime(S)/(x_0))$ is free as a $\mathbb{Z}_2[[t]]$-module, generated by the generalized Kauffman states defined in \cite[Figure~4]{ozsszcube}, and concentrated in a single algebraic grading. The unreduced $H_\ast(C^\prime(S))$ is also concentrated in a single algebraic grading. To calculate the structure of $H_\ast(C^\prime(S))$ as an $\cR^\prime$-module, we use a planar Heegaard diagram for $S$ defined exactly as in \cite[Figure~9]{ozsszcube} with extra bivalent vertices of the layered diagram treated as if they had come from smoothing a crossing.  So, the diagram looks like that on the right in Figure~\ref{fig:hgddiagram} near any bivalent vertex. The chain complex coming from the planar diagram is well-defined over $\mathbb{Z}_2[\underline{x}(S),t]$ (no completions required) because the planar diagram satisfies a stronger admissibility property than the modified initial diagram. However, we consider it over the larger ring $\cR^\prime$ because we need to compare its homology to the homology of the chain complex coming from the modified initial diagram. The same procedure of handleslides and destabilizations described in the proof of \cite[Lemma~3.7]{ozsszcube} shows that the two chain complexes are quasi-isomorphic. The planar diagram has a canonical generator, which is a cycle, defined by making the same choice of intersection point near each vertex as Ozsv\'ath and Szab\'o do in \cite[Proposition~3.10]{ozsszcube}. Incoming differentials from chains with algebraic grading one higher than the canonical generator produce all of the quadratic local relations, the linear local relations associated to bivalent vertices, and the non-local relations that appear in the definition of $\cA(S)$. Since $H_\ast(C^\prime(S))$ is concentrated in a single algebraic grading, this completes the calculation and establishes the isomorphisms claimed in~(\ref{singcalc}) and~(\ref{singcalcreduced}).

Now consider a layered braid diagram $D$ with $m$ crossings, and let $D_I$ denote its $I$-resolution. The spectral sequence constructed in \cite[Section~4]{ozsszcube} comes from a filtration on $$V(D)=\bigoplus_{I\in\{0,1\}^m}H_\ast\!\left(\underline{CFK}^-(D_I)\otimes M_I\right),$$ where $M_I$ is the Koszul complex on linear relations coming from all vertices in diagram $D_I$.  To define the filtration, Ozsv\'ath and Szab\'o consider a planar Heegaard diagram that simultaneously encodes each possible state (positive, negative, singularized, smoothed) of a crossing \cite[Figure 12]{ozsszcube}.  To adapt this Heegaard diagram to $D$, we need only add a small piece like that shown on the right in Figure~\ref{fig:hgddiagram} near any bivalent vertex. Call the diagram from \cite[Figure~12]{ozsszcube} so adapted the \emph{master diagram}.

Using particular choices of generators near crossings in the master diagram, Ozsv\'ath and Szab\'o define a filtration on $V(D)$. They also define maps that count holomorphic disks intersecting certain regions near crossings in the master diagram \cite[Section~4]{ozsszcube}. In \cite[Proposition~5.2]{ozsszcube}, they verify that some of these maps (those with the appropriate gradings) are the same as the edge maps in Section~\ref{sec:diffl} of this paper, under the identification of $H_\ast(C^\prime(D_I))$ with $\cA_I(D)$. The description of all of the maps on $V(D)$ and the proof of \cite[Proposition~5.2]{ozsszcube} depend only on the form of their Heegaard diagram near crossings, so they apply unchanged to our master diagram.  Taken together, the maps defined by counting appropriate holomorphic disks near crossings in the master diagram form an endomorphism of $V(D)$. Lemma~4.6 of \cite{ozsszcube} shows that $V(D)$ with this endomorphism is quasi-isomorphic to the chain complex $\underline{CFK}^-(D)$, which is the twisted knot Floer homology of the classical knot $D$, defined via the traditional holomorphic disks construction and regarded as an $\mathbb{Z}_2[x_0][[t]]$-module.  Again, the arguments depend only on the properties of the master diagram near crossings in $D$, so they carry through unchanged to our situation. Therefore, as in \cite[Theorem~4.4]{ozsszcube}, the filtration on $V(D)$ gives rise to a spectral sequence with $E_1$ page $$\bigoplus_{I\in\{0,1\}^m}H_\ast\!\left(\underline{CFK}^-(D_I)\otimes M_I\right),$$ with $d_1$ differential the zip and unzip maps defined algebraically, and converging to $\underline{HFK}^-(D)$.

Finally, in Section~5, Ozsv\'ath and Szab\'o argue that this spectral sequence collapses after the $E_1$ stage for grading reasons. The gradings in this paper are defined identically to those in \cite{ozsszcube}, so the same argument shows that the spectral sequence here collapses.  The immediate result is an isomorphism of $\mathbb{Z}_2[x_0][[t]]$-modules
$$H_\ast\left(\bigoplus_{I\in\{0,1\}^m}H_\ast\!\left(\underline{CFK}^-(D_I)\right)\otimes M_I\right)\isom H_\ast\!\left(\underline{CFK}^-(D)\right)$$
Inverting $t$ in the ground ring throughout the spectral sequence, then applying the isomorphism from~(\ref{singcalc}) allows us to identify the left side with the cube of resolutions complex $C(D)$ used in this paper:
$$H_\ast\left(C(D)\otimes_{\cR[\underline{x}(D)]}\widehat{\cR[\underline{x}(D)]}\otimes\mathbb{Z}_2\right)\isom H_\ast\!\left(\underline{CFK}^-(D)\otimes_{\mathbb{Z}_2[[t]]}\mathbb{Z}_2[t^{-1},t]]\right)$$
A standard theorem about twisted coefficients in knot Floer homology, stated as \cite[Lemma~2.2]{ozsszcube}, completes the identification with $H_\ast\left(CFK^-(D)\right)\otimes_{\mathbb{Z}_2}\mathbb{Z}_2[t^{-1},t]]$. The reduced statement follows similarly.
\end{proof}

\bibliographystyle{abbrv}
\bibliography{gilmoremathbib}

\end{document}